\documentclass[10pt]{amsart}
\usepackage{amsmath}
\usepackage{latexsym}
\usepackage{amssymb,amsthm,amsfonts}
\usepackage{graphicx}
\usepackage[active]{srcltx}
\usepackage{color} 
\usepackage{esint}

\renewcommand{\a }{\alpha }
\renewcommand{\b }{\beta }
\renewcommand{\d }{\delta }

\newcommand{\D }{\Delta }

\newcommand{\e }{\varepsilon }
\newcommand{\g }{\gamma}
\newcommand{\G }{\Gamma }
\renewcommand{\l }{\lambda }
\renewcommand{\L }{\Lambda }

\newcommand{\n }{\nabla }
\newcommand{\var }{\varphi }

\newcommand{\s }{\sigma }
\newcommand{\Sig }{\Sigma}
\renewcommand{\t }{\tau }
\renewcommand{\th }{\theta }

\renewcommand{\O }{\Omega }

\newcommand{\ov}{\overline}
\newcommand{\intbar}{\mathop{\int\makebox(-13.5,0){\rule[4pt]{.7em}{0.3pt}}%
\kern-6pt}\nolimits}

\newcommand{\be}{\begin{equation}}
\newcommand{\ee}{\end{equation}}
\newenvironment{pf}{\noindent{\sc Proof}.\enspace}{\rule{2mm}{2mm}\medskip}
\newenvironment{pfn}{\noindent{\sc Proof}}{\rule{2mm}{2mm}\medskip}

\newcommand{\R}{\mathbb{R}}
\newcommand{\C}{\mathbb{C}}

\newcommand{\Z}{\mathbb{Z}}

\newcommand{\N}{\mathbb{N}}
\newcommand{\pa}{\partial}
\newcommand{\dis}{\displaystyle}

\newcommand{\al}{\alpha}

\newcommand{\scp}{\scriptscriptstyle}
\renewcommand{\ov}[1]{\overline{#1}}
\newcommand{\un}[1]{\underline{#1}}
\renewcommand{\i}{\infty}
\newcommand{\sg }{\sigma }
\newcommand{\eps}{\varepsilon}
\newcommand{\sd}{\mathbb{S}^2}
\newcommand{\dt }{\delta }
\newcommand{\graf}[1]{\left\{\begin{array}{ll}#1\end{array}\right.}

\usepackage{accents}
\newcommand{\beq}{\begin{equation}}
\newcommand{\eeq}{\end{equation}}
\newcommand{\fr}{\frac }
\newcommand{\abs}[1]{\left\vert#1\right\vert}

\newtheorem{lem}{Lemma}[section]
\newtheorem{pro}[lem]{Proposition}
\newtheorem{thm}[lem]{Theorem}
\newtheorem{rem}[lem]{Remark}
\newtheorem{cor}[lem]{Corollary}

\begin{document}

\title[An improved geometric inequality] {An improved geometric inequality via vanishing moments, with applications to singular Liouville equations}

\author{Daniele Bartolucci and Andrea Malchiodi}

\address{Universit\`a di Roma {\em Tor Vergata}, Dipartimento di Matematica,
Via della Ricerca Scientifica n. 1 00133 Rome, Italy  and SISSA, via Bonomea 265,
34136 Trieste, Italy}

\email{bartoluc@axp.mat.uniroma2.it, malchiod@sissa.it}



\begin{abstract}
We consider a class of singular Liouville equations on compact surfaces motivated
by the study of Electroweak and Self-Dual Chern-Simons theories, the Gaussian curvature
prescription with conical singularities and Onsager's description of turbulence.
We analyse the problem of existence variationally, and show how the angular distribution
of the conformal volume near the singularities may lead to improvements in the Moser-Trudinger inequality, and in turn to lower bounds on the Euler-Lagrange functional.
We then discuss existence and non-existence results.
\end{abstract}

\maketitle

\section{Introduction}\label{s:in}

\noindent On a compact orientable surface $(\Sig, g)$ without boundary
and with metric $g$ we consider the equation
\begin{equation}\label{eq:e-1}
    - \Delta_g u = \rho \, \left( \frac{h(x)e^{2u}}{\int_{\Sig}
        h(x)e^{2u} d V_g} - a(x) \right)
    - 2 \pi \sum_{j=1}^m \alpha_j  \left(\delta_{p_j} - \frac{1}{|\Sig|}\right), \qquad
    \quad \int_{\Sig} a(x) dV_g = 1.
\end{equation}
Here $\rho$ is a positive parameter, $a, h : \Sig \to \R$ two smooth
functions, $h(x) > 0$ for every $x \in \Sigma$, $\alpha_j > 0$, $p_j \in \Sig$ and $|\Sig|$ denotes the area of
$\Sig$, that is $|\Sig|=\int\limits_{\Sig}d V_g$.\\

The analysis of \eqref{eq:e-1} is motivated by the study of vortex type configurations
(see for example  \cite{bardem}, \cite{btcmp02}, \cite{cl3}, \cite{sy2} and \cite{cl4}, \cite{djlw2},
\cite{linweiye}, \cite{maru2}, \cite{noltarcmp}, \cite{osuz}, \cite{st}) in the Electroweak theory of Glashow-Salam-Weinberg \cite{lai},
and in Self-Dual Chern-Simons theories \cite{dunne}, \cite{hkp}, \cite{jawe}. We refer the reader to
\cite{cl25}, \cite{tardcds} and to the monographs \cite{tar}, \cite{yang} for further details and an up to date set of references
concerning these applications.
Other classical problems call up for the study of \eqref{eq:e-1} such as the prescribed  Gaussian curvature problem
on surfaces with conical singularities \cite{bdm}, \cite{cl}, \cite{Troy0}, \cite{tro} and  Onsager's
statistical mechanics description of two-dimensional turbulence \cite{clmp2} in presence of vortex sources
\cite{chakie}. Moreover the
study of \eqref{eq:e-1} and of the corresponding Dirichlet problem
(see \eqref{prob:Diri} below) on bounded domains $\O\subset \R^2$ has an independent interest related to
the description of rotational shear flows \cite{turya} and/or Euler flows in presence of vortex sources \cite{barmon-1}.
It turns out that the structure of the solutions' set for these elliptic problems is extremely sensitive to the data.
For example it is well known (see Proposition \ref{p:nonexdisk} below) that if $\O$ is the unit ball with just
one singularity (i.e. $m=1$) located at the origin and $h\equiv 1$,
then we have non existence of solutions for \eqref{prob:Diri} with $\rho\geq 4\pi(1+\al_1)$.
On the other side, if either $\Sig\equiv \sd$ or $\O$ is simply connected
(at least to our knowledge, with few exceptions which will be shortly discussed below), there are no general existence results at all for either \eqref{eq:e-1} or \eqref{prob:Diri}
with $\rho>4\pi$. In any case there are still relatively few results at hand concerning
existence of solution for \eqref{eq:e-1} (some of which will be described in more detail below).
The reason for this gap consists essentially in the well known issue of non coercivity of the variational
functionals (see \eqref{eq:I_rho} below) associated to the study of these problems when $\rho\geq 4\pi$.
It seems in particular that we have a quite unsatisfactory understanding of some of
these existence/non existence problems as for example a model case \cite{barmon-1} suggests that the Dirichlet problem on simply connected domains
should admit at least one solution for each $\rho \in (0, 4 \pi \min_{i = 1, \dots, m} \{1+\alpha_i\}) \setminus 4 \pi \N$.
It is one of our motivations to fill this gap here.

\

\noindent In \cite{btcmp02} an existence theorem is
proved via min-max methods for surfaces with positive genus and for $\rho \in (4\pi, 8\pi)$.
More recently the latter result has been extended in \cite{bardem}, \cite{bdm}, still for positive
genus, for $\rho$ outside the discrete {\em critical} set (see \eqref{eq:quantsing} and Theorem \ref{t:bt2} below)
found in \cite{btcmp02}. In \cite{maru} the case of arbitrary genus was treated, but only for $\a_j \leq 1$ for all $j$
and $\rho \in (4\pi, 8\pi)$.

In \cite{cl4} existence results are deduced by calculating the Leray-Schauder degree when  $\alpha_j \geq 1$ for all $j$ and $\rho \in (4\pi,8\pi)$. In \cite{cl25}, \cite{cl3} an on-going project to compute the  Leray-Schauder degree of the equation is presented, using refined blow-up analysis and Lyapunov-Schmidt reductions, concerning the
case $\a_j \in \N$ for all $j$. This approach has been successful for the
{\em regular case}, when all $\a_j$'s are zero: a formula for the degree of the equation
has been derived in \cite{cl1}, \cite{cl2} building upon previous
blow-up analysis and quantization results in \cite{bm}, \cite{li}, \cite{lisha}.

Necessary and sufficient conditions for the existence of a solution for the Dirichlet problem
in the critical case $\rho=4\pi$ in presence of singularities on simply connected domains has been recently found in
\cite{barlin-2} (see also \cite{barlin}).

We also refer to \cite{dpem}, \cite{esp} for some perturbative results providing solutions of multi-bump type (via implicit function theorems) for special values of the parameter $\rho$ for \eqref{eq:e-1} on bounded
two dimensional domains with Dirichlet boundary conditions.

\

\noindent The goal of this paper is to prove a key inequality for treating \eqref{eq:e-1} variationally in general situations, and to present some applications to the existence problem
in simple cases. In particular
{\em 1)} compared to \cite{bardem}, \cite{bdm} and \cite{btcmp02} we remove the assumption on the
genus; {\em 2)} compared to \cite{cl4} and \cite{maru} we also allow the $\a_j$'s and $\rho$ to be arbitrarily large;
{\em 3)} compared to \cite{cl25}, \cite{cl3} we do not require $\al_j \in \N$: we notice that the structure of solutions to \eqref{eq:e-1} might
change drastically depending on $\rho$ when the coefficients $\al_j$ are not integer, see Remark \ref{r:commenti} {\em (a)}, and that in some situations one might still have existence
when the degree of the equation vanishes (see \cite{dj});  {\em 4)} compared to \cite{dpem}, \cite{esp} and \cite{barlin-2}, we allow generic values of $\rho$.
We also expect that our result, combined with those in \cite{cama}, might allow to
treat the case of $\al$'s with arbitrary sign, more relevant for geometric applications, as well as give precise homological
information on the variational structure (about the latter issue see more comments after
\eqref{eq:quantsing}).

\

\noindent Problem \eqref{eq:e-1} admits an equivalent formulation with variational
structure: letting $G_p(x)$ denote the Green's function of $- \D_g$ on
$\Sig$ with pole at $p$, i.e. the unique solution to
\begin{equation}\label{eq:Green}
   - \D_g G_p(x) = \d_p - \frac{1}{|\Sig|} \quad \hbox{ on }
   \Sig, \qquad \quad  \mbox{ with } \quad \int_{\Sigma} G_p(x) \, dV_g =0,
\end{equation}
by the substitution
\begin{equation}\label{eq:change}
    u \mapsto u + 2\pi \sum_{j=1}^m \alpha_j G_{p_j}(x),\qquad h(x)\mapsto \tilde{h}(x)= h(x) e^{-2\pi \sum_{j=1}^m \alpha_j G_{p_j}(x)}
\end{equation}
\eqref{eq:e-1} transforms into an equation of the type
\begin{equation}\label{eq:ee}
    - \Delta_g u  = \rho \left (\frac{\tilde{h}(x) e^{2u}}{\int_\Sig
    \tilde{h}(x) e^{2u} dV_g} - \frac{1}{|\Sig|} \right ) \qquad \quad \hbox{ on } \Sig.
\end{equation}
In general the constant $\frac{1}{|\Sig|}$ in \eqref{eq:ee} is replaced by a
smooth function $\tilde{a}(x)$ with $\int_\Sig \tilde{a}(x) dV_g = 1$: this term is indeed
rather harmless, and we will not comment on this issue any further.

Since $G_p$ has the asymptotic behavior $G_{p_j}(x) \simeq
\frac{1}{2\pi} \log \frac{1}{d(x,p_j)}$ near $p_j$, by \eqref{eq:change}
the function $\tilde{h}$ satisfies
\begin{equation}\label{eq:tildeh}
    \tilde{h} > 0 \hbox{ on } \Sig \setminus \cup_j \{p_j\}; \qquad \quad
  \qquad \tilde{h}(x) \simeq \g_j d(x,p_j)^{2 \alpha_j}  \quad \hbox{ near }
  p_j
\end{equation}
for some constant $\g_j > 0$, where $d(\cdot, \cdot)$ stands for the distance 
induced by $g$.

\

\noindent Problem \eqref{eq:ee} is the Euler-Lagrange equation of the
functional
\begin{equation}\label{eq:I_rho}
    I_{\rho,\underline{\a}}(u) = \int_{\Sig} |\n_g u|^2 dV_g + 2 \frac{\rho}{|\Sig|} \int_\Sig  u \,
    dV_g - \rho \log \int_\Sig \tilde{h}(x) e^{2u} dV_g; \qquad u \in H^{1}(\Sig).
\end{equation}
One basic tool for treating such kind of functionals is the well known
Moser-Trudinger inequality
\begin{equation}\label{eq:mtin}
    \log \int_\Sig e^{2(u - \ov{u})} dV_g \leq \frac{1}{4 \pi}
    \int_\Sig |\n_g u|^2 dV_g + C; \quad \qquad  u \in H^1(\Sig),
    \quad  \ov{u} = \fint_\Sig u \, dV_g,
\end{equation}
see e.g. \cite{moser} and \cite{Fon}. The value $\frac{1}{4 \pi}$ is sharp
in \eqref{eq:mtin}, as one can  insert in the above inequality a test function like
\begin{equation}\label{eq:bubble}
    \var_{\l,x}(y) = \log \frac{\l}{1 + \l^2 dist(x,y)^2}; \qquad \quad
    \l > 0, x \in \Sigma,
\end{equation}
and check that both sides diverge to infinity at the same rate.
This function is usually called a {\em standard bubble},
since the conformal metric $\tilde{g} = e^{2 \var_{\l ,x}}$ endows $\Sig$
with a spherical metric near $x$.

In presence of singularities, namely when a weight $\tilde{h}$ as in
\eqref{eq:tildeh} multiplies the exponential term, a modified sharp
Moser-Trudinger inequality was derived in \cite{chenwx} and \cite{tro}
(see also \cite{cygc}, or also \cite{cia}, \cite{fm} for exensions to general 
settings), and takes the form
\begin{equation}\label{eq:t}
  \log \int_\Sig \tilde{h} e^{2(u-\ov{u})} dV_g \leq \frac{1}{4 \pi}
      \frac{1}{\min \left\{ 1, \min_j \{1 +\a_j\} \right\}}
     \int_\Sig |\n u|^2 dV_g + C.
\end{equation}
As one can see, the constant is bigger when one of the coefficients - say
$\a_{\ov{j}}$ - is negative, as $\tilde{h}$ is singular near $\a_{\ov{j}}$.
However, when all of the $\a_j$'s are positive, as in the case we are considering,
the best constant remains $\frac{1}{4\pi}$. The fact that $\tilde{h}$ is equal
to zero
at all singular points does not give a smaller constant, as one may initially
expect: inserting in the inequality a {\em bubble} at a regular point $x$
does not pick up any effect of the vanishing of $\tilde{h}$ near the $p_j$'s.

\

\noindent From \eqref{eq:t} one has that $I_{\rho,\underline{\a}}$ is bounded from below
for $\rho < 4 \pi$, and hence one can find solutions of \eqref{eq:e-1} by
globally minimizing $I_{\rho,\underline{\a}}$, which is coercive, using the direct methods
of the calculus of variations. When $\rho > 4 \pi$ instead the situation
becomes more delicate, as $\inf I_{\rho,\underline{\a}} = - \infty$: one might however
try to obtain solutions as saddle points. We describe next some previous results in the literature which rely on this strategy.

Even though $I_{\rho,\underline{\a}}(u)$ is not bounded below on $H^1(\Sig)$,
one might hope to find suitable conditions on $u$ to recover some control. Calling
$\mathfrak{A} \subseteq H^1(\Sig)$ a set of functions for which this lower bound holds,
one can then try to show that $\mathfrak{A}$ is always
intersected along a suitable family of min-max maps.

For the regular case of \eqref{eq:e-1} such a lower bound was obtained in
by W.Chen and C.Li in \cite{cl} (extending previous results in
\cite{aubin} and \cite{mo} for the standard sphere) under the condition that two
subsets of $\Sigma$ with positive mutual distance both contain a finite
portion of the total mass. Under such an assumption, one finds that
the best constant in \eqref{eq:mtin} can be chosen arbitrarily close to
$\frac{1}{8 \pi}$: as a consequence, when $\rho < 8 \pi$ and when
$I_{\rho,\underline{\al}}(u)$
is large negative, $e^{2u}$ has to concentrate near a single point of $\Sig$
(similarly, if $\rho < 4 (k+1) \pi$, $e^{2u}$ concentrates near at most
$k$ points, as shown in \cite{dm}). This property was used in \cite{djlw} to obtain existence on
surfaces of positive genus when $\rho \in (4 \pi, 8 \pi)$, and in
\cite{dj} (relying on an argument in \cite{dm} for the $Q$-curvature prescription
problem) for $\rho \not\in 4 \pi \N$ on all surfaces. The restriction $\rho \not \in
4 \pi \N$ is a compactness condition which allows to apply the deformation lemma,
see \cite{str}, \cite{lucia} (see also \cite{cl2}, \cite{cl3} for some results
concerning the case $\rho \in 4 \pi \N$).

For the singular case, a related approach has been used in \cite{btcmp02}
where, through a new quantization property (see Theorem \ref{t:bt2}) the result
in \cite{djlw} was extended to the case of positive $\a_j$'s (and, still, for positive
genus). In particular, compactness is obtained provided $\rho \not\in \L_{\underline{\al}}$,
where
\begin{equation}\label{eq:quantsing}
  \L_{\underline{\al}} = \left\{ 4 k \pi + 4\pi\sum_{j \in J} (1 + \a_j) \; | \;
  k \in \N \cup \{0\}, J \subseteq \{1, \dots, m\} \right\}.
\end{equation}
The latter existence result was later generalized in \cite{bardem}, \cite{bdm} to the case of
arbitrarily large (but positive) values of $\rho$.

\

\noindent The last two existence results however do not fully capture the variational
features of the presence of the singularities, from three different aspects.
They do neither extend to the case of the sphere or to the case when some
negative weights are present (which could be relevant, we recall, for the Gaussian curvature
prescription and to the study of turbulent flows interacting with vortex sinks). Finally, the analysis is not sufficient if one wants to fully characterize
from the homological point of view the structure of sublevels of $I_{\rho, \underline{\al}}$, to
compute for example the degree of the equation as it has been done in \cite{mald}
for the regular case.

An improvement of inequality \eqref{eq:t}, more intrinsically related to the presence
of singularities, was derived in \cite{det}: it was shown that for any $\a >
-1$ there exists $C_\a > 0$ such that
  \begin{equation}\label{eq:det}
  \log \int_B |x|^{2 \a} e^{2(u - \ov{u})} dV_g \leq \frac{1}{4 (1+\a) \pi}
    \int_B |\n_g u|^2 dV_g + C_\a; \qquad \quad u \in H^1_r(B),
  \end{equation}
where $B$ is the unit ball of $\R^2$ and $H^1_r$  the space of radial
functions in $B$ of class $H^1$. This result has a previous related counterpart in
\cite{mc}, where the case of curvatures with $\Z_k$ symmetry and polynomial decay
in $\R^2$ was considered, among others.

In \cite{maru} a general improvement (without assuming any symmetry) was found 
for $\a \in (0,1]$
and $\rho \in (4 \pi, 8 \pi)$: recall that in this case, by the above discussion,
a low energy for $u$ implies concentration of the volume near at most one point.
The novelty in \cite{maru} was to derive an extra characterization of this
point, which takes into account both the {\em scale} of concentration of the
volume measure and its {\em center of mass}. More precisely, it was proven that
there exists a continuous map $\beta$ from low sublevels of $I_{\rho, \underline{\al}}$
into $B$ such that if $\beta(u)$ hits the singularity then \eqref{eq:t} holds
with $\frac{1}{4 \pi}$ replaced by $\frac{1+\e}{4 \pi (1+\a)}$, where
$\e$ can be chosen arbitrarily small (see Proposition \ref{l:mr} for
more details). Notice that this condition relaxes the radiality in \cite{det} to a
two-dimensional constraint, and that it allows an arbitrarily small
scale of concentration at a single point (so \cite{cl} would not apply). The
condition $\a \leq 1$ is indeed sharp in this argument, as one can find
counterexamples for $\a > 1$. We also refer to \cite{maru2} where a somehow
related strategy is used for {\em Toda systems} (arising from non-abelian
theories).

\

\noindent The main goal of this paper is to find a general condition to get an improved
inequality for arbitrary $\al$'s, with no symmetry requirements, and which is
flexible enough to be combined with min-max arguments. As we will try to describe,
our approach combines the scaling invariance properties in \cite{maru} and
the possibility of volume concentration at multiple points
(as in \cite{dj}, \cite{dm}).

To explain this condition in more detail, suppose we are on the
unit ball $B$ of $\R^2$ and that we are dealing with only one singularity at the
origin with weight $\al$. Let $\tilde{f}_u$ denote the probability
measure on $B$
  \begin{equation}\label{eq:tfu}
  \tilde{f}_u = \frac{\tilde{h}(x) e^{2u}}{\int_B \tilde{h}(x) e^{2u} dx}.
  \end{equation}
Roughly speaking, our result can be interpreted as a version of the above concentration property at finitely-many points in a \underline{complete} setting,
blowing-up the  metric near the singularity as $\tilde{g} = \frac{1}{|x|^2} (dx)^2$ so that the Euclidean metric becomes cylindrical. To state this property rigorously, assuming that
$\rho \in (4 k \pi, 4 (k+1) \pi)$, given $\d > 0$ small we define
    \begin{equation}\label{eq:intk}
    J_{k,\d}(\tilde{f}_u) = \sup_{x_1, \dots, x_k \neq 0} \int_{\cup_{i=1}^k
    B_{\d |x_i|}(x_i)} \tilde{f}_u dx.
    \end{equation}
To describe our strategy, we first consider two alternatives which may occur:
when  $J_{k,\d}(\tilde{f}_u)$ is close to $1$ and when it is not.

When this quantity is close to $1$, we are in a situation similar
to Chen and Li's (but in the cylindrical metric). For the regular case,  the argument in \cite{dj} (or in \cite{dm}) implies that if $k$ small balls in $B$ contain most
of the volume (as in this first alternative), then it is possible to find a continuous map from these measures into the {\em formal barycenters of $B$ of
order} $k$, namely the probability measures of the form
$B_k = \left\{ \sum_{i=1}^k t_i \d_{x_i} \; : \; t_i \geq 0,
\sum_{i=1}^k t_i = 1, x_i \in B \right\}$.

In our case, it is natural to incorporate the dilation invariance of
the problem (corresponding to a translation along the axis of the
cylinder), and to project onto the barycenters of order $k$ of $S^1$,
which coincides with the cylinder factoring out the translations.

For doing this, we define the probability measure on the circle
 \begin{equation}\label{eq:mumuuu}
  \mu_u(A) = \int_{\tilde{A}} \tilde{f}_u dx;  \qquad \quad
  A \subseteq S^1, \quad \tilde{A} = \cup_{t \in (0,1]} t A.
 \end{equation}
When $J_{k,\d}(f)$ is close to $1$ then $\mu_u \simeq \sum_{i=1}^k t_i
\d_{\th_i}$ (in the distributional sense) for some $t_i \geq 0$ and some
$\th_i \in S^1$. The $k$-barycenters of $S^1$, $(S^1)_k$, are known to be homotopically equivalent to
$S^{2k-1}$, see Theorem 1.1 and Corollary 1.5 in \cite{kk}. It is however convenient for us
to understand this set in more detail, proving that it is indeed {\em homeomorphic}
to a (piecewise smooth) immersed sphere $\mathcal{S}_k$ in $\C^k$ with interior $\mathcal{U}_k$ being a  neighborhood of the origin, see Section \ref{s:proj}.
This is useful in order to construct a continuous projection of a small neighborhood $\mathcal{N}_k$ of $\mathcal{S}_k$ onto $\mathcal{S}_k$ itself.

More precisely, let
  \begin{equation}\label{eq:FkFk}
  F_k(\tilde{f}_u) = \left( \int_{S^1} z \, d \mu_u, \int_{S^1} z^2 d \mu_u, \dots,
  \int_{S^1} z^k d \mu_u \right),
  \end{equation}
mapping the probability measures on $S^1$ into $\C^k$.
Using this map, we define $\mathcal{S}_k$ to be $F_k((S^1)_k)$, which
can be seen to be a homeomorphism, and we check that in the first
alternative ($J_{k,\d}(\tilde{f}_u)$ close to $1$) $F_k(\tilde{f}_u)$ lies in
$\mathcal{N}_k$, so we can project continuously onto
$\mathcal{S}_k \simeq (S^1)_k$.

\

\noindent We consider next the second of the two alternatives, namely when
$J_{k,\d}(\tilde{f}_u)$ is
bounded away from $1$. One thing to be immediately noticed
is that while in a compact situation one always obtains
weak convergence of a sequence of probability measures to a
probability measure, in the complete case some part of the {\em mass}
(or the whole one) might disappear by {\em vanishing}.

We show by a covering argument  that if $J_{k,\d}(\tilde{f}_u)$ is
bounded away from $1$ then either some volume concentrates near at
least $k+1$ well separated points with respect to the cylindrical metric, or 
that some part of the
measure spreads on the cylinder (giving rise to a vanishing), see Lemma
\ref{l:altprel}. In either case, using harmonic liftings
and some argument in \cite{maru} (see Proposition \ref{l:mr}) we show
that the constant in \eqref{eq:t} improves by a factor $\min
\{1+k, 1+\a\}$, see Proposition \ref{p:impr} (recall that $\al$ stands for the
weight of the singularity at the origin).
One condition, easy to verify, which guarantees this improvement is the vanishing of the moments of the measure $\mu_u$ up to order $k$, see Corollary \ref{c:moments}. Qualitatively,
this is quite similar to  the requirement on $\beta(u)$ in \cite{maru},
see the comments after \eqref{eq:det}. Furthermore, it would apply to a symmetric
case as in \cite{mc}, but it only imposes finitely-many integral constrains on $u$.

\

\noindent We employ the previously described  inequality to find new existence
results for \eqref{eq:e-1} and its analogue Dirichlet problem on bonded domains $\O\subset\R^2$, that is
\begin{equation}\label{prob:Diri}
 \graf{
 -\Delta u=\rho \dfrac{\tilde{h} e^{2u}}{\int_{B} \tilde{h} e^{2u}dx} & \mbox{in}\quad \O\\
 u=0 & \mbox{on}\quad \pa \O,
}
\end{equation}
where, $\tilde{h}(x)= h(x) e^{-2\pi \sum_{j=1}^m \alpha_j G_{p_j,0}(x)}$ for some strictly positive and smooth $h$ on
$\O$ and, for $p\in\O$, $G_{p,0}(x)$ denotes the Green's function uniquely defined by
\begin{equation}\label{eq:Green:Diri}
   - \D G_{p,0}(x) = \d_p\quad \hbox{ on }
   \O, \qquad \quad  G_{p,0}(x)=0\quad\hbox{ on }\;\pa \O.
\end{equation}

\smallskip

\begin{thm}\label{t:exdisk}
Let $\O \subseteq \R^2$, and let $m \neq 0$. Then problem \eqref{prob:Diri}
admits a solution for every $\rho \in (0, 4 \pi \min_{i = 1, \dots, m} \{1+\alpha_i\}) \setminus 4 \pi \N$.
\end{thm}

\smallskip

\begin{thm}\label{t:sphere} Suppose $\Sig$ is a topological sphere, and
let $m \geq 2$. Then \eqref{eq:e-1} has a solution provided $\rho \in
(0, 4 \pi \min_{i = 1, \dots, m} \{1+\alpha_i\}) \setminus 4 \pi \N$.
\end{thm}

\smallskip

\begin{rem}\label{r:commenti} (a) The upper bounds on $\rho$ in Theorems
\ref{t:exdisk} and \ref{t:sphere} are sharp:
in Section \ref{s:ex} we complement our results with Propositions \ref{p:nonexdisk}
and \ref{p:nonexsph}, giving non existence for larger values of $\rho$ (for
$m = 1$ in the unit ball or $\R^2$ and for $m = 2$ on $S^2$).

(b) We will prove in detail Theorem \ref{t:exdisk} only for simply
connected domains and for $m = 1$, since this case is the simplest
one requiring our new estimates. We sketch the argument for the other
cases in Remark \ref{r:ab}, since the proof adapts quite easily. The counterpart 
of Theorem \ref{t:sphere} for surfaces with positive genus (actually a 
more general version of it, without upper bounds on $\rho$)
was proved in \cite{bdm}: in Remark \ref{r:ab} we will briefly discuss how
our method can also be adapted to these surfaces when $\rho < 4 \pi \min_{i = 1, \dots, m} \{1+\alpha_i\}$.
\end{rem}

\smallskip

\noindent To prove Theorem \ref{t:exdisk} we employ a min-max scheme which uses
the formal barycenters of $S^1$. More precisely, let $k$ be the unique integer
for which $\rho \in (4 k \pi, 4 (k+1) \pi)$. Then, since $J_{k,\d}(\tilde{f}_u)$
separate from $1$ leads to a lower bound for $I_{\rho, \underline{\al}}$, the
above discussion suggests that if $u$ has low energy, then the measure $\mu_u$
on $S^1$ (see \eqref{eq:mumuuu}) should be close to some element of $(S^1)_k$ in the distributional sense.

We build a min-max scheme based on this idea: starting from $\s \in (S^1)_k$
we define a test function (see \eqref{eq:testdisk}) for which the associated conformal
volume resembles $\s$, and on which $I_{\rho, \underline{\al}}$ attains large negative
values, see Lemma \ref{l:lowen}. Since $(S^1)_k \simeq S^{2k-1}$ (see the comments before
\eqref{eq:FkFk}), the family of these test functions forms
a $2k-1$-dimensional sphere in $H^1_0(B)$, on which the supremum of $I_{\rho, \underline{\al}}$
is very low.

We then complete this family with a map from a topological ball in $\C^k$ into
$H^1_0(B)$, and we consider the min-max value
associated to this construction, see Proposition \ref{p:crit0}. The improved inequality in Proposition \ref{p:impr}
is used to show that the min-max value is strictly larger than the maximal value
at the boundary, otherwise by Proposition \ref{p:proj} we would be able to find
(naively) a retraction from the unit ball of $\C^k$ onto its boundary, which is
a contradiction. Details
are given in Section \ref{s:ex}.

The compactness issue due to a lack of knowledge about
the Palais-Smale condition can be tackled via by now standard means: varying the parameter
$\rho$ and reasoning as in \cite{lucia}, \cite{str} it is possible to find
$\rho_n \to \rho$ for which $I_{\rho_n, \underline{\al}}$ has a bounded Palais-Smale
sequence and hence a solution. Convergence can then be proved using Theorem \ref{t:bt2}
and Corollary \ref{c:comp}.

The proof of Theorem \ref{t:sphere} can be handled using minor modifications:
the main point is that concentration of conformal volume may occur either near the
singularity $p_1$ or near $p_2$. The
corresponding improved inequality is given in Proposition \ref{p:imprsph}, where one can
see that both weights $\al_1, \al_2$ play a role. The proof of the non existence results
in Propositions \ref{p:nonexdisk} and \ref{p:nonexsph} are shown using well-known
Poho\u{z}aev type identities (see also the recent paper \cite{carl} for non existence
results on surfaces with positive genus).

\

\noindent The paper is organized as follows. In Section
\ref{s:pr} we list some preliminary results on elementary inequalities,
the Moser-Trudinger inequality and some of its improvements, together
with the compactness result from \cite{btcmp02}. In Section \ref{s:proj}
we show how to embed continuously the barycenters $(S^1)_k$ into $\C^k$
using moments of measures on the unit circle, and how to project continuously onto this image the
family of functions for which $J_{k,\d}$ is close to $1$.
In Section 4 we then analyse the complementary situation, proving a
dichotomy result in Proposition \ref{l:altprel} and then the
improved inequality in both alternatives. In Section \ref{s:ex} we then
prove both existence and non existence results.

 \

\begin{center}

 {\bf Acknowledgements}

\end{center}

\noindent The authors are supported by the project FIRB-Ideas {\em
Analysis and Beyond} and by the MiUR project {\em Variational methods and
nonlinear PDEs}. The authors are grateful to G.Tarantello for the interpretation
of the improved inequality in Proposition \ref{p:impr} in terms of
vanishing momenta., to D.Ruiz for pointing out to us Appendix E in \cite{bredon}
and to A.Carlotto for some useful comments.

\section{Notation and preliminaries}\label{s:pr}

\noindent This section contains some useful preliminary material,
including elementary inequalities, some
variants of the Moser-Trudinger inequality from \cite{cl}, \cite{dj}, \cite{dm},
\cite{maru}, and a compactness results from \cite{btcmp02}.

We will deal with either compact Riemannian surfaces $(\Sig, g)$, with
or without boundary, or with the unit ball $B$ of $\R^2$. We let $d(x,y)$
be the distance of two points $x, y$, while $B_r(p)$ will stand for
the open metric ball of radius $r$ and center $p$. We also set, for convenience
$$
B_r=\{x\in B \,|\, d_g(x,0)<r\}, \qquad \hbox{ and } \qquad
A(s,t)=\{x\in B \,|\, s<d_g(x,0)<t\}.
$$
The symbol
$\fint_\O u \, dV_g$ denotes the average integral $\frac{1}{|\O|}
\int_\O u \, dV_g$.  For $\a > 0$ we set
\begin{equation}\label{eq:ka}
  k_\a = \min \left\{ k \in \N \; | \; k \geq \a \right\}.
\end{equation}
If $u \in H^1(\Sig)$ or $u \in H^1_0(B)$, and if $\O$ has smooth boundary and
is compactly contained in the domain of $u$, we denote by
$\mathcal{H}_\O(u)$ the harmonic extension of $u$ inside $\O$, namely
we set
\begin{equation}\label{eq:hlift}
 \mathcal{H}_\O (u) = \begin{cases}
   v & \hbox{ in } \O ; \\
   u & \hbox{ in } \Sig \setminus \O \quad (\hbox{resp. } B \setminus \O),
   \end{cases}
\end{equation}
where $v$ is the solution of
$$
   \begin{cases}
   \D _g v = 0 & \hbox{ in } \O; \\
   v = u & \hbox{ on } \partial \O.   \end{cases}
$$
For two probability measures $\mu_1, \mu_2$ defined on $B$
the {\em Kantorovich-Rubinstein} distance is defined as
\begin{equation}\label{eq:dflat}
  d_{KR}(\mu_1, \mu_2) = \sup_{\|f\|_{Lip} \leq 1} \left|
  \int_B f \, d \mu_1 - \int_B f  \, d \mu_2 \right|.
\end{equation}
If $\tilde{h}$ is as in \eqref{eq:tildeh}, for $u$ of class $H^1$ we will set
\begin{equation}\label{eq:fu}
   f_u = \tilde{h}(x) e^{2u}; \qquad \qquad \tilde{f}_u =
   \frac{\tilde{h}(x) e^{2u}}{\int_B \tilde{h}(x) e^{2u} dx}.
\end{equation}
We will use  similar notations for functions which are defined on a
compact surface $\Sig$.

\

\noindent Generic large positive constants are always denoted by $C$, $\tilde{C}$,
etc.: even though we
allow constants to vary, we will often stress their dependence on other
constants or parameters, as we need sometimes to be careful in the ordering of their
choices.

\

\subsection{Some elementary inequalities}\label{ss:el}

We begin with two elementary Lemmas: the first can be proved e.g. using Fourier analysis, while the second follows from Poincar\'e's inequality.

\begin{lem}\label{l:ext} Let $p \in B$, let $s > 0$, and suppose $B_{2s}(p)
\subseteq B$. For  $u \in H^1(B_{2s}(p))$, let $\mathcal{H}_{B_s(p)}(u)$
be as in \eqref{eq:hlift}. Then there exists a universal constant $C_0$,
independent of $u$, $p$ and $s$, such that
$$
  \int_{B_s(p)} |\n \mathcal{H}_{B_s(p)}(u)|^2 dx \leq C_0 \int_{B_{2s}(p)
  \setminus B_s(p)} |\n u|^2 dx.
$$
\end{lem}

\

\begin{rem} \label{r:bsbordo}
A similar result holds when $u \in H^1_0(B)$, $B_s(p) \cap \partial B
\neq \emptyset$ and $d(p, \partial B) < \frac{s}{2}$, with $s \leq \frac{1}{4}$:
this condition controls from below the angle formed by $\pa B_s(p)$ and $\pa B$
at their intersection points.
Analogous statements can be also proved for small metric balls on a given compact
surface.
\end{rem}

\

\begin{lem}\label{l:medie} Let $p \in B$, $s > 0$, and suppose $B_s(p)
\subseteq B$. Let $\ov{C}_0$ be a fixed constant, and suppose that
$B_r(q) \subseteq B_s(p)$ with $r \geq \ov{C}_0^{-1} s$, and
$d(B_r(q), \partial B_s(p)) \geq \ov{C}_0^{-1} s$.
Let $u \in H^1(B_s(p))$:  then there exists another constant
$\tilde{C}_0$, depending on $\ov{C}_0$ but independent of $s$, $q$, $r$ and $u$, such that
$$
  \left| \fint_{\partial B_r(q)} u \, d \s - \fint_{\partial B_s(p)}
  u\,  d \s \right| \leq \tilde{C}_0 \|\n u\|_{L^2(B_s(p))}.
$$
\end{lem}

\

\subsection{Improved Moser-Trudinger inequalities}

We start by recalling the well known Moser-Trudinger inequality for
surfaces with or without boundary (see e.g. \cite{cygc}, \cite{Fon}, \cite{moser}).

\begin{pro} Let $\Sig$ be a compact surface. Then
\noindent \begin{enumerate} \item[a)] If $\Sig$ has no boundary,
\begin{equation}\label{eq:mtdisk}
    \log \int_{\Sig} e^{2 u} dV_g \leq \frac{1}{2 \pi} \int_{\Sig}
  |\n_g u|^2 dV_g + 2 \fint_{\Sig} u \, dV_g + C \quad \hbox{ for every } u \in
  H^1(\Sig).
\end{equation}
\item[b)] If $\Sig$ has boundary,
\begin{equation}\label{eq:mtdisk0}
    \log \int_{\Sig} e^{2 u} \ dV_g \leq \frac{1}{4 \pi} \int_{\Sig}
  |\n_g u|^2 dV_g + C \quad \hbox{ for every } u \in H^1_0(\Sig).
\end{equation}
\end{enumerate}
\end{pro}

\noindent As we remarked in the introduction, the constant $\frac{1}{4\pi}$ in
\eqref{eq:mtdisk} is sharp, as one can see plugging in the inequality
test functions as in \eqref{eq:bubble}.

The next result, proven in  \cite{cl} for $\ell = 1$ and in \cite{dj},
\cite{dm} for general $\ell$, gives a criterion
for getting a smaller multiplicative constant in the Moser-Trudinger
inequality. Basically, it asserts that the more $e^{2u}$ is {\em spread},
the smaller constant can be chosen in \eqref{eq:mtdisk}. The proof,
not reported here, relies on localizing the Moser-Trudinger inequality
(through suitable cut-off functions) near the sets, called $\Omega_i$'s,
on which there is concentration of volume.

\begin{pro}\label{l:imprc0} Let $\Sig$ be a compact surface with no boundary, $\tilde{h}:
\Sigma \to \R$ with $0 \leq \tilde{h}(x) \leq C_0$. Let $\Omega_1, \dots,
 \Omega_{\ell+1}$ be subsets of $\Sig$ with $dist(\Omega_i,\Omega_j)
\geq \delta_0$ for some $\delta_0>0$ if $i \neq j$, and fix $\g_0 \in \left( 0,
\frac{1}{\ell+1} \right)$. Then, for any $\e > 0$ there exists a
constant $C = C(C_0, \e, \delta_0, \g_0, \ell)$ such that
$$
  \log \int_\Sig f_u dV_g \leq C + \frac{1+\e}{
  4 (\ell + 1) \pi} \int_{\Sig} |\n_g u|^2 dV_g + 2 \fint_{\Sigma} u \,
  dV_g
$$
for all functions $u \in H^1(\Sig)$ satisfying
\begin{equation}\label{eq:ddmmi}
    \int_{\Omega_i} \tilde{f}_u dV_g
    \geq \g_0, \qquad \quad \quad  i =1,\dots, \ell+1.
\end{equation}
A similar result holds, without the average of $u$ on the right-hand side,
if $\Sig$ has boundary and $u \in H^1_0(\Sig)$.
\end{pro}

\noindent The case of surfaces with boundary is not explicitly written in
\cite{cl} but their proof can be adapted with minor changes to cover this situation as well.
A useful corollary of Proposition \ref{l:imprc0} is the following, which describes
the set of functions for which the Euler-Lagrange functional is large negative.
For the proof, which uses Proposition \ref{l:imprc0} and a covering
argument, see \cite{dj} and \cite{dm} (see also \cite{cl} or \cite{djlw} for $k = 1$).

\begin{cor}\label{c:conc} Suppose $\rho < 4 (k+1) \pi$. Then, given any
$\e, r > 0$ there exists $L = L(\e,r) > 0$ such that
$$
  I_{\rho, \underline{\a}}(u) \leq - L \qquad \Rightarrow \qquad \dis
  \int_{\cup_{j=1}^k B_{r}(x_j)} \tilde{f}_u \,
  dV_g > 1 - \e \quad
  \hbox{ for some } x_1, \dots, x_k \in \Sig.
$$
\end{cor}

\

\noindent The next improved Moser-Trudinger inequality is established in \cite{maru},
and exploits the role played by the singularities in a more subtle way. While Proposition
\ref{l:imprc0} is based  on the separation of concentration regions, Proposition \ref{l:mr} involves a separation in the {\em scales} of concentration.

\begin{pro}\label{l:mr} Consider the case of one singularity at the origin in $B$
(in our previous notation, $m = 1$ and $p_1 = 0$), and let $\a = \a_1$.
Let $\eta$ be a small positive constant, and fix $\t > 0$.
Let $u \in H^1_0(B)$, and suppose  there exists $s \in \left(0, \frac{1}{4} \right)$ such that
$$
   \int_{s < |x| < 4 s} |\n u|^2 dx \leq \eta \int_B |\n u|^2 dx,
$$
and such that
  \begin{equation}\label{eq:ininin}
  \int_{|x| < s} \tilde{f}_u dx \geq \t; \qquad \qquad \int_{|x| > 4 s} \tilde{f}_u
  dx \geq \t.
  \end{equation}
Then, there exists a universal constant $C_0 > 0$ and $\tilde{C} =
\tilde{C}(\eta, \t, \a)$ (independent of $s$) such that one has the inequality
$$
 (1 + \a) \log \int_B f_u dx \leq \frac{1}{4 \pi} \left( \a \int_{|x| < 2 s} |\n u|^2 dx +
 \int_{|x| > 2 s} |\n u|^2 dx + C_0 \eta \int_B |\n u|^2 dx\right) + \tilde{C}.
$$
\end{pro}

\begin{pf} The details of the proof can be found in Proposition 4.1 of \cite{maru}:
for the reader's convenience, since we will also need some modified version of this
result (see Remark \ref{r:centro}) we will sketch here the main arguments.

First of all, we modify $u$ in $B_{4s} \setminus B_{s}$ so it becomes constant
 in $B_{3s} \setminus B_{2s}$: precisely, if we let
$$
  \chi_s(r) = \min \left\{ \frac{1}{s} (r-s), 1, \frac{1}{s}
  (4s-r) \right\}; \qquad \quad \hat{u}(s) = \fint_{B_{4s}
  \setminus B_s} u \, dx,
$$
and define
$$
  \tilde{u}(x) =
  \begin{cases}
  \chi_s(|x|) \hat{u}_s +
    \left( 1 - \chi_s(|x|) \right) u(x) & \hbox{ for } x \in B_{4s}
      \setminus B_s; \\
      u(x) & \hbox{ for } x \in B \setminus (B_{4s}     \setminus B_s).
  \end{cases}
$$
By Poincar\'e's inequality and our assumptions we have that (choosing possibly
a larger universal $C_0$)
\begin{equation}\label{eq:gradtildeu2}
  \int_{B_{4s} \setminus B_s} |\n \tilde{u}|^2 dx \leq C_0 \eta \int_B |\n u|^2 dx.
\end{equation}
Hence, by the first inequality in \eqref{eq:ininin}, the asymptotics of $\tilde{h}$,
a change of variables (a dilation bringing $B_{2s}$ into $B$), by
\eqref{eq:mtdisk0} (used on $\tilde{u} - \hat{u}(s)$), and \eqref{eq:gradtildeu2}
one finds
\begin{equation}\label{eq:mmmm}
   \log \int_{B} f_u dx \leq \frac{1}{4 \pi} \left( \int_{B_s}
   |\n u|^2 dx + C_0 \eta \int_B |\n u|^2 dx \right) + 2 \hat{u}(s)
   + 2 (1 + \a) \log s + \tilde{C}.
\end{equation}
Moreover one has
\begin{equation}\label{eq:stupid}
    \int_B f_u dx \leq \frac 1 \t \int_{B \setminus B_{2s}}f_{\tilde{u}} dx =
    \int_{B \setminus B_{2s}} \frac{\tilde{h}}{|x|^{4\a}} |x|^{4 \a} e^{2 \tilde{u}}
    dx \leq \frac{C}{s^{2\a}} \int_{B \setminus B_{2s}} e^{2v} dx,
\end{equation}
with $v(x) = \tilde{u}(x) + 2 \a w(x)$, where
$$w(x)= \left \{ \begin{array}{ll} \log (2 s) & x \in B_{2s},
\\  \log |x| & x \in B \setminus B_{2s}. 
 \end{array} \right.
$$
Notice that
\begin{eqnarray*}
  \int_{B \setminus B_{2s}} |\n v|^2 dx & = &
  \int_{B \setminus B_{2s}} |\n \tilde{u}|^2 dx + 4 \a^2
  \int_{B \setminus B_{2s}} \frac{1}{|x|^2} dx \\
  & + & 4 \a \int_{B \setminus B_{2s}} \langle \n \tilde{u}, \n (\log|x|)
  \rangle dx.
\end{eqnarray*}
We integrate by parts to obtain
\begin{equation}\label{fund}
\int_{B \setminus B_{2s}} |\n v|^2 dx \leq \int_{B \setminus B_{2s}}
|\n \tilde{u}|^2 dx + 8 \pi \a^2 \log \frac 1 s -
8 \pi \a \fint_{\partial B_{2s}} \tilde{u} \, d\s + \tilde{C}.
\end{equation}
Next, by using the second inequality in \eqref{eq:ininin}, \eqref{eq:stupid},
\eqref{eq:mtdisk0} for $v$ and the fact that $\tilde{u} \equiv \hat{u}(s)$
on $\pa B_{2s}$, we get
\begin{equation}\label{eq:mmmm2222}
\log \int_B f_u dx \leq 2 \a \log \frac 1 s + \frac{1}{4\pi}
\int_{B \setminus B_{2s}} |\n \tilde{u}|^2 dx + 2 \a^2
\log \frac 1 s - 2 \a \hat{u}(s) + C_0 \eta
\int_{B} |\n u|^2 dx + \tilde{C}.
\end{equation}
By using \eqref{eq:mmmm} together with \eqref{eq:mmmm2222} we finally deduce
$$
  (1 + \a) \log \int_B f_u dx \leq \frac{\a}{4\pi} \int_{B_s}
  |\n u|^2 dx + \frac{1}{4\pi} \int_{B \setminus B_{3s}} |\n u|^2
  dx + C_0 \eta \int_{B} |\n u|^2 dx + \tilde{C},
$$
which is the desired conclusion.
\end{pf}

\begin{rem}\label{r:centro} (a) The above proposition also works when the
center of the ball is shifted by an amount of order of the
radius. More precisely, if we have the same assumptions of Proposition
\ref{l:mr} replacing $B_s$ (resp. $B_{4s}$) by $B_s(p)$ (resp. $B_{4s}(p)$)
with $|p| \leq \ov{C} s$, the same result will hold provided $B_{4s}(p)
\subseteq B$, and allowing the dependence of $\tilde{C}$ also on $\ov{C}$.
The proof follows the same lines as before (combined in particular with Lemma
\ref{l:medie}) and requires minor adaptations
from \cite{maru}, where this case is treated on compact surfaces.

(b) The same assertion as in the previous part of this remark holds if the
condition $B_{4s}(p) \subseteq B$ is replaced by $d(p, \pa B) \leq \frac{s}{2}$
(see Remark \ref{r:bsbordo}). To see this, one can simply use a localization
argument for \eqref{eq:mtdisk0} as in the proof in \cite{cl}, see the
comments before Proposition \ref{l:imprc0}.
\end{rem}

\subsection{Compactness of solutions}

Concerning \eqref{eq:ee}, we have the following result, proved in
\cite{btcmp02} via blow-up analysis, extending previous theorems in
\cite{bls}, \cite{bm} and \cite{lisha} for the regular case (see also
\cite{barmon} for the case of negative $\a_i$'s).

\begin{thm}\label{t:bt2}
Let $\Sig$ be a compact surface, and let $u_i$ solve \eqref{eq:ee}
with $\tilde{h}$ as in \eqref{eq:tildeh}, $\rho = \rho_i$, $\rho_i
\to \ov{\rho}$, with $\alpha_j
> 0$ and $p_j \in \Sig$. Suppose that $\int_\Sig
f_{u_i} dV_g \leq \ov{C}$ for some fixed $\ov{C}
> 0$. Then along a subsequence $u_{i_k}$ one of the following
alternative holds:
\begin{description}
  \item[(i)] $u_{i_k}$ is uniformly bounded from above on $\Sig$;
  \item[(ii)] $\max_\Sig \left( 2 u_{i_k} - \log \int_\Sig f_{u_{i_k}}
 dV_g  \right) \to + \infty$ and there exists a finite blow-up set
   $S = \{ q_1, \dots, q_l\} \in \Sig$ such that

   $(a)$ for any $s \in \{1, \dots, l\}$ there exist $x^s_k \to
   q_s$ such that $u_{i_k}(x^s_k) \to + \infty$ and $u_{i_k} \to
   - \infty$ uniformly on the compact sets of $\Sig \setminus
   S$,

   $(b)$ $\rho_{i_k} \tilde{f}_{u_{i_k}} \rightharpoonup \sum_{s=1}^l
   \b_s \delta_{q_s}$ in the sense of measures, with $\b_s = 4
   \pi$ for $q_s \neq \{p_1, \dots, p_m\}$, or $\b_s = 4 \pi (1
   + \alpha_j)$ if $q_s = p_j$ for some $j = \{1, \dots, m\}$.
   In particular one has that
   $$
  \ov{\rho} = 4 \pi n + 4 \pi \sum_{j \in J} (1 + \alpha_j),
   $$
for some $n \in \N \cup 0$ and $J \subseteq \{1, \dots, m\}$
(possibly empty) satisfying $n + |J| > 0$, where $|J|$ is the
cardinality of the set $J$.
\end{description}
A similar result holds on bounded domains, for Dirichlet boundary data.
\end{thm}

\noindent From the above result we obtain immediately the
following corollary, which will be useful to prove our existence results.
Recall the definition of $k_\a$ in \eqref{eq:ka} and of $\Lambda_{\underline{\al}}$
in \eqref{eq:quantsing}.

\begin{cor}\label{c:comp}
Consider problem \eqref{eq:ee} in $B$, with Dirichlet boundary data, and
suppose $m = 1$. Let $\rho \in (4\pi, 4 (k_\a + 1) \pi)$. Then the set of
solutions is uniformly bounded in $C^2(\ov{B})$ provided
$$
  \rho \neq 4 \pi (1 + \alpha) \qquad   \hbox{ and }
  \qquad \rho \neq 4 k \pi, \quad   k = 1, \dots, k_\a.
$$
\end{cor}

\

\section{Momenta of probability measures and the set of formal barycenters of $S^1$}\label{s:proj}

\noindent For fixed $k\in\N$ we denote by $\un{z}_k\in\C^k$ the vector $\un{z}_k=(z_1,\ldots,z_k)$, by $\mathcal{D}_k$
the following subset of $\C^k$
$$
\mathcal{D}_k =\{\un{z}_k\in\C^k\,|\,|z_1|=|z_2|=\cdots=|z_{\scp k}|=1\},
$$
and for $R>0$
$$
B_R^{(k)}=\{\un{z}_k\in\C^k\,|\,|z_1|^2+|z_2|^2+\cdots+|z_{\scp k}|^2<R^2\}.
$$
We also set
$$
\R^{(+)}_k:=\{\un{t}_k\in \R^k\,:\,t_i> 0,\;\forall\,i=1,\cdots,k\},
$$
and
$$
S_k:=\left\{\un{t}_k\in [0,1]^k\,:\,\sum\limits_{i=1}^{k}t_i=1\right\},\quad
\accentset{\circ}{S}_k:=\left\{\un{t}_k\in (0,1)^k\,:\,\sum\limits_{i=1}^{k}t_i=1\right\}.
$$
For $\un{t}_k\in S_k$ we denote by $\sg_k$ an element in the space of $k$-baricenters of $S^1$, that is
$$
(S^{1})_k\ni\sg_k = \sum_{i=1}^k t_i \d_{\th_i},\quad \th_i\in[0,2\pi),\,\forall\, i=1,\cdots,k,
$$
and finally define $F_k:(S^1)_k\mapsto \C^k$ to be the following map,
$$
  F_k(\sg_k) = \left( \int_{S^1} z \, d\sg_k, \int_{S^1} z^2 d \sg_k, \dots,
  \int_{S^1} z^k d \sg_k   \right).
$$

\begin{pro}\label{p:proj} There exist small constants $\d_k, \t_k > 0$ with the following property.
If $\d \leq \d_k $ there exists a continuous (with respect to the Kantorovich-Rubinstein distance) map $\Xi_k$ from the set of functions $f \in L^1$, $\int_B f dx = 1$
satisfying
\begin{equation}\label{11-05-12.1}
    J_{k,\d}(f) > 1 - \t_k,
\end{equation}
into $\mathcal{S}_k$. Moreover if $f_n \to \s \in (S^1)_k$ in the Kantorovich-Rubinstein metric, then
$\Xi_k(f_n) \to F_k(\s)$.
\end{pro}

\noindent The proof of Proposition \ref{p:proj} can be deduced as a direct consequence of the following,

\begin{pro}\label{p:homeo} The map $F_k$ realizes a homeomorphism (with respect to the Kantorovich-Rubinstein metric)
between $(S^1)_k$ and a topological sphere $\mathcal{S}_k$ in $\C^k$ which bounds a neighborhood $\mathcal{U}_k$
of $0 \in \C^k$.
\end{pro}

\noindent We first use Proposition \ref{p:homeo} to prove Proposition \ref{p:proj}.

\

\begin{pfn} {\sc of Proposition \ref{p:proj}} It is straightforward to check
that for $\d_k, \t_k > 0$ small enough
and for any $\d\leq \d_k$ then any $f$ satisfying \eqref{11-05-12.1} is close
(with respect to the Kantorovich-Rubinstein metric) to a $k$-barycenter of $S^1$, and hence it
is mapped in some neighborhood $\mathcal{N}_k$ of $\mathcal{S}_k$.

By Theorem E.3 in \cite{bredon}, since $\mathcal{S}_k$ is a topological sphere,
it is a retract of some neighborhood of its in $\C^k$. Choosing
$\d_k$ and $\t_k$ possibly smaller, we find that $\mathcal{N}_k$ will
be contained in this neighborhood.
The map $\Xi_k$ is finally obtained as the composition of $F_k$ with the above retraction.
\end{pfn}

\

\begin{pfn} {\sc of Proposition \ref{p:homeo}}
For $\un{t}_k\in \R^{(+)}_k\cup S_k$  let $\Psi_{\scp k,\un{t}_k}:
\C^k\mapsto \C^k$ be defined by
$$
\Psi_{\scp k,\un{t}_k}(\un{z}_k)=\left(
\begin{array}{c}
t_1 z_1+t_2 z_2+\ldots+ t_k z_k\\
t_1 z_1^2+t_2 z_2^2+\ldots+ t_k z_k^2\\
\vdots\\
t_1z_1^k+t_2 z_2^k+\ldots+t_k z_k^k
\end{array}
\right).
$$
Hence,
$$
F_k\left(\left(S^1\right)_k\right)=\left\{\Psi_{\scp
k,\un{t}_k}\left(\mathcal{D}_k\right)\right\}_{\un{t}_k\in S_k}.
$$
Let $\deg\left( \Psi, B^{(k)}_R,\un{y}_k\right)$ denote the topological degree (see for example \cite{Zeidler})
of a map $\Psi:\C^k\mapsto\C^k$ relative to $B^{(k)}_R$ with
respect to $\un{y}_k\in \C^k$. We have the following:

\begin{lem}\label{bardem-lem} (see Lemma 3.1 in \cite{bardem})
For fixed $\un{t}_k\in \R^{(+)}_k$ there holds,
\beq\label{deg-k-bdem} \deg\left( \Psi_{\scp k,\un{t}_k},
B^{(k)}_R,\un{0}_k\right)=k!. \eeq
\end{lem}

\noindent Although Lemma 3.1 in \cite{bardem} concerns the case
$R=1$ the proof provided there works indeed for general $R$.

\

\noindent Let us consider the new variables $w_i=t_iz_i\in
B^{(1)}_1,\;\forall\,i\in\{1,\ldots,k\}$, so that for
$\un{t}_k\in\accentset{\circ}{S}_k$ and $\un{z}_k\in
\mathcal{D}_k$ we have in particular $t_i=|w_i|\in (0,1)$ and
$\un{w}_k\in B^{(k)}_1$. Hence $\Psi_{\scp k,\un{t}_k}$ takes the
form
$$
\Phi_k(\un{w}_k):=\Psi_{\scp k,\un{t}_k}(\un{z}_k)=\left(
\begin{array}{ccccccc}
w_1&+& w_2&+&\ldots &+&w_k\\
& \vdots  & & \vdots &  & \vdots &\\
|w_1|^{(1-j)}w_1^j&+&|w_2|^{(1-j)} w_2^j&+&\ldots &+& |w_k|^{(1-j)} w_k^j\\
& \vdots  & & \vdots &  & \vdots &\\
|w_1|^{(1-k)}w_1^k &+& |w_2|^{(1-k)}w_2^k &+& \ldots &+&
|w_k|^{(1-k)} w_k^k
\end{array}
\right)
$$
and we conclude that
$$
F_k\left(\left(S^1\right)_k\right)=\left\{ \Psi_{\scp
k,\un{t}_k}\left(\mathcal{D}_k\right)\right\}_{\un{t}_k\in S_k}
\equiv \Phi_k(\pa \mathcal{R}_k),
$$
where
$$
\mathcal{R}_k=\{\un{w}_k\in\C^k\,|\,|w_1|+|w_2|+\cdots+|w_{\scp
k}|<1\}.
$$
We compute next a topological degree related to \eqref{deg-k-bdem}.
\begin{lem}\label{l:deg-k}
We have \beq\label{deg-k} \deg\left( \Phi_{k},
\mathcal{R}_k,\un{0}_k\right)=k!, \eeq and in particular
$\Phi_{k}(\un{w}_k)=0\,\iff\,\un{w}_k=\un{0}_k$.
\end{lem}

\begin{pfn}
By using Lemma 3.2 in \cite{bardem} we see that it is enough to
prove the assertion with $\mathcal{R}_k$ replaced by $B^{(k)}_1$, that is
$$
\deg\left( \Phi_{k}, B^{(k)}_1,\un{0}_k\right)=k!.
$$
In view of \eqref{deg-k-bdem}, for $s\in[0,1]$ we set
\begin{footnotesize}\begin{align*}
\mathcal{H}_k(\un{0}_k,0):=&\un{0}_k,\hfill\\
\mathcal{H}_k(\un{w}_k,s):=&
\begin{scriptsize}\fr1{sk+(1-s)}\left(
 \!\!\begin{array}{ccccccc}
w_1&\!\!+&\!\!w_2&\!+&\!\cdots&\!\!+&\!\! w_{\scp k}\\
&\scp{\vdots}& &\scp{\vdots}& &\scp{\vdots}&\\
w_1\left(\fr{w_1}{\sqrt{s+(1-s)\abs{w_1}}}\right)^{j-1}&\!\!+&\!\!w_2\left(\fr{w_2}{\sqrt{s+(1-s)\abs{w_2}}}\right)^{j-1}&+&\cdots&\!\!+&\!\! w_{\scp k}\left(\fr{w_{\scp k}}{\sqrt{s+(1-s)\abs{w_{\scp k}}}}\right)^{j-1}\\
&\scp{\vdots}& &\scp{\vdots}& &\scp{\vdots}&\\
w_1\left(\fr{w_1}{\sqrt{s+(1-s)\abs{w_1}}}\right)^{k-1}&\!+&\!w_2\left(\fr{w_2}{\sqrt{s+(1-s)\abs{w_2}}}\right)^{k-1}&+&\cdots&\!+&\! w_{\scp k}\left(\fr{w_{\scp k}}{\sqrt{s+(1-s)\abs{w_{\scp k}}^2}}\right)^{k-1}\\
\end{array}
 \right)\!\!,\end{scriptsize}
\end{align*}\end{footnotesize}
the latter being defined for $\un{w}_k\neq\un{0}_k$. Clearly
$\mathcal{H}_k$ is continuous and
$$
\mathcal{H}_k(\un{w}_k,1)=\Psi_k(\un{w}_k):=\Psi_{\scp
k,\un{t}_k}(\un{w}_k)\left.\right|_{\un{t}_k=\un{\frac{1}{k}}_k},\quad
\mathcal{H}_k(\un{w}_k,0)=\Phi_k(\un{w}_k).
$$
Therefore, to obtain \eqref{deg-k}, we are left to prove that
$\mathcal{H}_k(\un{w}_k,s)\neq\un{0}_k$ for any $s\in[0,1]$ and
for any $\un{w}_k\in\pa B^{(k)}_1$. More generally, we will show
by induction that
$$
\qquad\textrm{if for some $\un{w}_k\in\C^k$ and  $s\in[0,1]$ it
holds $H_k(\un{w}_k,s)=\un{0}_k$, then
$\un{w}_k=\un{0}_k$.}\qquad\quad(P)_{\scp k}
$$
We notice that for $s=1$ this is Lemma 3.2 in \cite{bardem}.
$(P)_1$ is trivially true, being $\mathcal{H}_1(\un{w}_1,s)=w_1$.
Now, suppose that for some $k\geq 2$ $(P)_m$ holds for any $m\leq
k-1$: then if by contradiction there exist $s\in[0,1]$ and
$\un{w}^{(1)}_k\in \C^k$ such that
$\mathcal{H}_k(\un{w}^{(1)}_k,s)=\un{0}_k$, we would have that
$w^{(1)}_i \neq 0$ for all $i\in\{1,\ldots,k\}$. Indeed otherwise,
if for some $i$ it holds $w^{(1)}_i=0$ then
$\mathcal{H}_{k-1}(\un{w}^{(1)}_{k-1,i})=\un{0}_k$, where
$\un{w}^{(1)}_{k-1,i}=(w^{(1)}_1,\ldots,w^{(1)}_{i-1},w^{(1)}_{i+1},\ldots,w^{(1)}_{\scp
k})$ is the $(k-1)$-tuple not including $w^{(1)}_i$. Hence
$(P)_{\scp k-1}$ implies $\un{w}^{(1)}_{k-1,i}=\un{0}_{k-1}$ and
we would conclude that
$\un{w}^{(1)}_k=\un{0}_k$ which is of course a contradiction.

At this point we are allowed to define: \beq\label{assurdo}
\tilde{w}^{(1)}_i:=\fr{w^{(1)}_i}{\sqrt{s+(1-s)\abs{w^{(1)}_i}}}\qquad\mbox{and}\qquad
t_i:=\fr{\sqrt{s+(1-s)\abs{w^{(1)}_i}}}{s k+(1-s)},\quad
i=1,\ldots,k. \eeq and conclude from the previous considerations
that
$\un{\tilde{w}}^{(1)}_k:=(\tilde{w}^{(1)}_1,\ldots,\tilde{w}^{(1)}_{\scp
k})\in\C^{k}\setminus\{\un{0}_k\}$. Of course
$\un{t}_k\in\R^{(+)}_k$ and since
$\mathcal{H}_k(\un{w}_k,s)=\un{0}_k$ if and only if
$\Psi_k(\un{\tilde{w}}^{(1)}_k)=\un{0}_k$ we deduce from Lemma 3.2
in \cite{bardem} that $\un{\tilde{w}}^{(1)}_k=\un{0}_k$ which is
the desired contradiction. \end{pfn}

\

\noindent Let $\Phi_{k}^{*}:\R^{2k}\mapsto \R^{2k}$ be the map $\Phi_k$ when
expressed in real coordinates and set
$$
\Upsilon_k=\{\un{w}_k\in\C_k\,|\,w_i=w_j,\;\;\mbox{for
some}\;i\neq j,\;\{i,j\}\subseteq\{1,\cdots,k\}\},
$$
$$
\Upsilon_k^{(0)}=\{\un{w}_k\in\C_k\,|\,w_j=0,\;\;\mbox{for
some}\;j\in\{1,\cdots,k\}\}.
$$
We have the following result
\begin{lem}\label{l:det} There holds
\begin{equation}\label{110512.1}
\mbox{\em det}\left(D\Phi_{k}^{*}(\mbox{\em
Re}(\un{w}_k),\,\mbox{\em Im}(\un{w}_k) )\right)\neq 0,\quad
\forall\un{w}_k\notin \Upsilon_k\cup \Upsilon_k^{(0)}.
\end{equation}
\end{lem}

\begin{pfn}
Setting $w_j=r_j e^{i \th_j}$, $\forall\,j=1,\cdots,k$, it is
straightforward to check that
$\mbox{det}\left(D\Phi_{k}^{*}\right)$ takes the form

$$
\mbox{det}\left(D\Phi_{k}^{*}(\mbox{Re}(\un{w}_k),\,\mbox{Im}(\un{w}_k))\right)=
r_1^k r_2^k\cdots r_k^k \mbox{det}\left(A_{2k}(\un{\th}_k)\right),
$$
where
$$
A_{2k}(\un{\th}_k) =\left(
\begin{array}{cccccccc}
\cos{(\th_1)}&\!\! \cos{(\th_2)}&\!\cdots&\!\!\cos{(\th_k)}&\!\!-\sin{(\th_1)}&\!\!-\sin{(\th_2)}&\!\cdots&\!\!-\sin{(\th_k)}\\
\sin{(\th_1)}&\!\! \sin{(\th_2)}&\!\cdots&\!\!\sin{(\th_k)}&\!\!\cos{(\th_1)}&\!\!\cos{(\th_2)}&\!\cdots&\!\!\cos{(\th_k)}\\
\cos{(2\th_1)}&\!\! \cos{(2\th_2)}&\!\cdots&\!\!\cos{(2\th_k)}&\!\!-2\sin{(2\th_1)}&\!\!-2\sin{(\th_2)}&\!\cdots&\!\!-2\sin{(\th_k)}\\
\sin{(2\th_1)}&\!\! \sin{(2\th_2)}&\!\cdots&\!\!\sin{(2\th_k)}&\!\!2\cos{(2\th_1)}&\!\!2\cos{(2\th_2)}&\!\cdots&\!\!\cos{(2\th_k)}\\
\vdots&\!\! \vdots&\!\!\vdots&\!\!\vdots&\!\cdots&\!\!\vdots&\!\! \vdots\\
\cos{(k\th_1)}&\!\! \cos{(k\th_2)}&\!\cdots&\!\!\cos{(k\th_k)}&\!\!-k\sin{(k\th_1)}&\!\!-k\sin{(k\th_2)}&\!\cdots&\!\!-k\sin{(k\th_k)}\\
\sin{(k\th_1)}&\!\! \sin{(k\th_2)}&\!\cdots&\!\!\sin{(k\th_k)}&\!\!k\cos{(k\th_1)}&\!\!k\cos{(k\th_2)}&\!\cdots&\!\!k\cos{(k\th_k)}\\
\end{array}
\right).
$$
It follows immediately that $\un{w}_k \notin \Upsilon_k^{(0)}$ is
a necessary condition for \eqref{110512.1}
to be satisfied.\\
Next observe that $\mbox{det}\left(A_{2k}(\un{\th}_k)\right)\neq
0$ whenever $\th_i\neq \th_j,\;\forall\,i\neq j$. In fact, letting
$\{\un{\mathrm{v}}_m\}_{m=1,\cdots,k}$ be the complex row vectors
$$
\un{\mathrm{v}}_m=\left(e^{m i \th_1},\cdots,e^{m i \th_k}\right),
$$
we see that $A_{2k}(\un{\th}_k)$ takes the form
$$
A_{2k}(\un{\th}_k) =\left(
\begin{array}{cc}
\un{\mathrm{v}}_1 & i \, \un{\mathrm{v}}_1 \\
\un{\mathrm{v}}_2 & 2 \, i \, \un{\mathrm{v}}_2 \\
\vdots\\
\un{\mathrm{v}}_k & k \, i \un{\mathrm{v}}_k
\end{array}
\right).
$$
Hence, if $\mbox{det}\left(A_{2k}(\un{\th}_k)\right)= 0$, then
there exists $\un{\l}_{k}\in\C^{k}$ such that,
$$
\sum\limits_{m=1}^{k}\l_m \left(
\begin{array}{c}
\un{\mathrm{v}}_m\\
m \, i \, \un{\mathrm{v}}_m
\end{array}
\right)=\un{0}_{2k},
$$
which readily implies that the complex vectors
$\{\un{\mathrm{v}}_m\}_{m=1,\cdots,k}$ must be linearly dependent.
Since the matrix whose rows are
$\{\un{\mathrm{v}}_m\}_{m=1,\cdots,k}$ is the restriction to
$\mathcal{D}_k$ of the Vandermonde-type matrix defined by
$\Psi_{k,\un{t}_k}(\un{z}_k)\left.\right|_{\un{t}_k=\un{1}_k}$,
where $\un{1}_k$ is the vector whose entries are all $1$, then a
well known argument shows that necessarily $\th_i=\th_j$ for some
$i\neq j$.
\end{pfn}

\

\noindent Since the region $\C^k\setminus \left( \Upsilon_k\cup\Upsilon_k^{(0)} \right)$ is
connected and
$\mbox{det}\left(D\Phi_{k}^{*}(\mbox{Re}(\un{w}_k),\,\mbox{Im}(\un{w}_k)
)\right)$ is continuous, it follows from Lemma \ref{l:det} that
$\mbox{det}\left(D\Phi_{k}^{*}\right)$ has in fact constant sign
(unless it is zero). In this situation, and by using
\eqref{deg-k}, one can conclude (see \cite{Zeidler} Vol. I p. 639
prob. 14.3d and also Theorem 14A and Corollary 14.8) that for
fixed $\un{b}_k\in \C^k$ the Vandermonde-type system
$$
\Phi_{k}(\un{w}_k)=\un{b}_k, \qquad\qquad\qquad (V)_{k}(\un{b}_k)
$$
admits at least one and at most $k!$ distinct solutions.

Hence $\Phi_k(\mathcal{R}_k)$ is open and we define
$\mathcal{U}_k:=\Phi_k(\mathcal{R}_k)$ and
$\mathcal{S}_k:=\pa\mathcal{U}_k$. The underlying idea toward the
conclusion of the proof is that $\sg_k \in
(S^1)_k\setminus(S^1)_{k-1}$ if and only if
$\un{t}_k\in\accentset{\circ}{S}_k$ and $\th_i\neq\th_j$
$\forall\,i\neq j$, that in terms of $\un{w}_k$ variables is
equivalent to
$\un{w}_k\in \pa \mathcal{R}_k\setminus \left( \Upsilon_k\cup\Upsilon_k^{(0)} \right)$.

\

\noindent For $i\in\{1,\ldots,k\}$ let $\Pi_{i}:\C^k\mapsto\C^k$ be the
standard permutation map
$$
\Pi_i(\un{w}_k)=\Pi_i((w_1,\ldots,w_i,w_{i+1},\ldots,w_k))=
\left(w_1,\ldots,w_{i+1},w_i,\ldots,w_k\right),
$$
which is defined with the periodic condition (for fixed $k\in\N$)
$k+1=1$. Clearly, for fixed $\un{b}_k\in \C^k$, each of the
$k!$ permutations of a given solution $\un{w}_k$ of
$(V)_{k}(\un{b}_k)$ will yield a distinct solution whenever
$\un{w}_k\notin \Upsilon_k$. If
$\un{b}_k\notin\Phi_k\left(\Upsilon_k\cup\Upsilon_k^{(0)}\right)$,
then to each of the corresponding $k!$ distinct solutions
$\un{w}_k^{(m)}\in\pa\mathcal{R}_k,\,\forall\,m=1,\cdots,k!,\;$
there correspond the same $k$-baricenter
$\sg_k=F_k^{-1}(\un{b}_k)=\sum\limits_{j=1}^k t_j \d_{\th_j}$,
where $t_j=|w_j^{(1)}|$, $\th_j=\mbox{arg}(w_j^{(1)})$. Clearly
$F_k$ is continuous, so Lemma \ref{l:det} and the Inverse Function
Theorem together imply that $F_k^{-1}$ is locally of class $C^1$
in $\Phi_k\left(\Upsilon_k\cup\Upsilon_k^{(0)}\right)$. We
conclude in particular that $F_k^{-1}$ is well defined and
continuous in $\mathcal{S}_k\setminus
\left( \Upsilon_k\cup\Upsilon_k^{(0)}\right)$. Next, we analyze
the case
$\un{a}_k\in\Phi_k\left(\pa\mathcal{R}_k\cap \left(\Upsilon_k\cup\Upsilon_k^{(0)}\right) \right)$.\\

\

\noindent {\bf Claim} If $\un{a}_k\in\Phi_k\left(\pa\mathcal{R}_k\cap
\left(\Upsilon_k\cup\Upsilon_k^{(0)}\right)\right)$ then to every solution
of $\Phi_k(\un{w}_k)=\un{a}_k$ there correspond a unique
$k$-baricenter $\s_{k}\in (S^1)_{k-1}\subset (S^1)_{k}$. In
particular $F_k^{-1}$ is well
defined and continuous on $\Phi_k\left(\pa\mathcal{R}_k\cap \left(\Upsilon_k\cup\Upsilon_k^{(0)}\right)\right)$.

\

\noindent \begin{pfn}
We argue by induction and observe that for $k=2$ it holds
$$
\Phi_2\left(\pa\mathcal{R}_2\cap\left(\Upsilon_2\cup\Upsilon_2^{(0)}\right)
\right)=\left\{
\left(\begin{array}{c}e^{i\th}\\e^{2i\th}\end{array}\right)
 \right\}_{\th\in[0,2\pi)}.
$$
It is readily seen that in fact to each $\un{a}_2\in
\Phi_2\left(\pa\mathcal{R}_2\cap\left( \Upsilon_2\cup\Upsilon_2^{(0)}\right)
\right)$
there correspond a unique $2$-baricenter $\s_2\in (S^1)_1\subset
(S^1)_2$. In particular $F_2^{-1}$ is continuous on $\left\{
\left(\begin{array}{c}e^{i\th}\\e^{2i\th}\end{array}\right)
 \right\}_{\th\in[0,2\pi)}$
with respect to the Kantorovich-Rubinstein metric since in this case we have
$$
F^{-1}\left({e^{i\th} \choose e^{2i\th}}\right)=\d_{\th}.
$$
Therefore, let us assume that the property in the statement of the
Claim holds for any $m\in\{1,\cdots,k-1\}$ and let us prove that
it holds for $m=k$ as well. Let $\un{w}_k\in
\pa\mathcal{R}_k\cap\left(\Upsilon_k\setminus \Upsilon_k^{(0)} \right)$
be such that $w_1=w_2$. Then set
$$
\tilde{w}_1=2w_1,\quad \tilde{w}_\ell=w_{\ell-1}, \quad \qquad \forall\,
\ell=3,\cdots,k.
$$
Hence it is well defined a $k-1$-dimensional vector
$\un{\tilde{w}}_{k-1}$ satisfying $\un{\tilde{w}}_{k-1}\in
\pa\mathcal{R}_{k-1}$. However it is not too difficult to verify
that any other $\un{w}_k$ satisfying either $\un{w}_k\in
\pa\mathcal{R}_k\cap \left( \Upsilon_k\setminus \Upsilon_k^{(0)}\right)$ and
one of the ${k\choose 2}$ constraints $w_i=w_j$ for some $i\neq j$
or $\un{w}_k\in
\pa\mathcal{R}_k\cap \Upsilon_k\cap \Upsilon_k^{(0)}$ can be
transformed in this way after an appropriate relabelling of the
indices. In particular a similar argument works for $\un{w}_k\in
\pa\mathcal{R}_k\cap \left( \Upsilon_k^{(0)}\setminus \Upsilon_k\right)$.
Hence

$$
\Phi_k\left(\pa\mathcal{R}_k\cap\left(\Upsilon_k\cup\Upsilon_k^{(0)}\right)\right)=
\left\{\Phi_k(\un{\tilde{w}}_{k-1})\right\}_{\un{\tilde{w}}_{k-1}\in\pa\mathcal{R}_{k-1}}.
$$
At this point let $\un{a}_k\in
\Phi_k\left(\pa\mathcal{R}_k\cap\left(\Upsilon_k\cup\Upsilon_k^{(0)}\right)\right)$
and define $\un{\tilde{a}}_{k-1}\in\C^{k-1}$ to be the vector
whose entries are the first $k-1$ entries of $\un{a}_k$ and
$\tilde{\Phi}_{k-1}(\un{\tilde{w}}_{k-1})$ the map whose rows are
the first $k-1$ rows of $\Phi_k$. By our discussion above and the
induction assumption, to any such $\un{\tilde{a}}_{k-1}$ there
correspond a unique $(k-1)$-baricenter
$\s_{k-1}=\s_{k-1}(\un{\tilde{a}}_{k-1})\in (S^1)_{k-1}$ which can
be obtained via any fixed solution of the system
$$
\tilde{\Phi}_{k-1}(\un{\tilde{w}}_{k-1})=\un{\tilde{a}}_{k-1}.
$$
In particular the inverse map determined in this way is
continuous. Finally, since $\un{\tilde{w}}_{k-1}\in
\pa\mathcal{R}_{k-1}$, it follows from Lemma \ref{l:deg-k} that
$\un{\tilde{a}}_{k-1}\neq \un{0}_{k-1}$. Hence, both the $k$-th
component of $\un{a}_k$ and the $(k-1)$-baricenter
$\s_k(\un{a}_k)\equiv \sg_{k-1}(\un{\tilde{a}}_{k-1})$ in the
pre-image of $\un{a}_k$ are fixed by $\un{\tilde{a}}_{k-1}$ and
therefore in particular the inverse map determined in this way is
continuous.
\end{pfn}

\bigskip

\noindent We can now conclude the proof of Proposition \ref{p:homeo}.
Putting $\un{a}_k\in \pa \mathcal{R}_k$, we set
$\Gamma_{s\un{a}_k}=\{s\un{a}_k\}_{s\in[0,1]}$ to be the ray
joining together the origin $\un{0}_k$ and $\un{a}_k$. We conclude
that if $\Gamma_{s\un{a}_k}\cap \Upsilon_k =\emptyset$ (which is
satisfied if and only if $\un{a}_k\notin \Upsilon_k$), then each
$\un{b}_k\in \Phi_k(\Gamma_{s\un{a}_k})$ admits exactly $k!$
distinct pre images and in particular that on any such
$\Gamma_{s\un{a}_k}$, $\Phi_k$ is injective. It is straightforward
to check by an induction argument as in the Claim that in fact
$\Phi_k$ is injective on $\Gamma_{s\un{a}_k}$ for $\un{a}_k\in
\Upsilon_k$ as well. Since $\Gamma_{s\un{a}_k}$ is a continuous
curve, we see that $\Phi_k(\mathcal{R}_k)$ is foliated by
$\left\{\Phi_k(\Gamma_{s\un{a}_k})\right\}_{\un{a}_k\in\pa\mathcal{R}_k}$.
Hence
$\mathcal{S}_k\equiv\left\{\Phi_k(\Gamma_{\un{a}_k})\right\}_{\un{a}_k\in\pa\mathcal{R}_k}
\equiv\Phi_k(\pa\mathcal{R}_k)=F_k((S^1)_k)$ is homeomorphic to a
$2k-1$ dimensional sphere embedded in $\C^k$.
\end{pfn}

\

\section{A general improved inequality}\label{s:im}

\noindent  Throughout this section we work on the unit ball $B$ and we  assume that
$$
 m = 1, \qquad p_1 = 0 \in B, \qquad \a = \a_1 > 0.
$$
The main result of this section is the following proposition,
which will be useful to obtain lower bounds on $I_{\rho, \underline{\a}}$ (see Corollary \ref{c:lowbd}).

\begin{pro}\label{p:impr} Let $\e > 0$, and let $k \in \{1, \dots, k_\a\}$. Suppose $\d_k$ and $\t_k$ are so small that Proposition
\ref{p:proj} applies, and let $\d \leq \d_k$.
Then there exists a constant $C_{\e,\a}$, depending only on $\e$ and $\a$, such that
$$
  \log \int_B f_u dx \leq \frac{1+\e}{4\pi \min\{ 1 + k, 1+\a\}}
  \int_B |\n u|^2 dx + C_{\e,\a}
$$
for all functions $u \in H^1_0(B)$ such that
$$
  J_{k,\d}(\tilde{f}_u) \leq 1 - \t_k.
$$
\end{pro}

\

\noindent The proof of the proposition is divided into several steps. We begin by
choosing a large constant $C_1$, depending on $\e$ and $\a $, such that
\begin{equation}\label{eq:C1}
    \frac{1}{\log C_1} = \frac{1}{32 (k_\a+1)^2}\frac{\e}{ 1 +  C_0^2},
\end{equation}
where $C_0$ is the constant in Lemma \ref{l:ext}. First,
we derive an alternative in case we are under the assumptions of Proposition
\ref{p:impr}. Consider the cylindrical metric as described in the Introduction,
after equation \eqref{eq:tfu}. Proposition \ref{l:altprel} asserts that if the
conformal volume is not concentrated near $k$ points of the cylinder obtained
from the blown-up metric, then either part of it accumulates near $k+1$ well separated
regions, or part of it {\em vanishes}. By this we mean that its integral
over bounded sets in some region of the cylinder becomes arbitrarily small.  The division
of the volume into $N$ parts in $(jj)$ is technical and will be needed in
the next subsections.

\begin{pro}\label{l:altprel}
Let $C_1$ be as in \eqref{eq:C1}, let  $k\in \{1, \dots, k_\a\}$ and let $f \in L^1(B)$
be such that $\int_B f \, dx = 1$ and $J_{k,\dt}(f)\leq 1-\tau_k$ ($\t_k$ and $\d$ are as in Proposition \ref{p:impr}). Then for any $\sg_0>0$ there exists $\sg\in (0,\sg_0]$,
depending on $\sg_0$, $\a$ and $\e$, but  not on $u$, such that the following
alternative holds: either

\

\noindent (j) there exist $k+1$ points $\{p_1,\cdots,p_{k+1}\}\subset B
\setminus \{0\}$ such that
\begin{equation}\label{eq:j1}
\int\limits_{B_{(10 C_1 k)^{-8}|p_i|}(p_i)} f \, dx \geq \s,\quad\forall\;i=1,\ldots,k+1,
\quad\mbox{and}\quad B_{(10 C_1 k)^{-4}|p_i|}(p_i)\cap B_{(10 C_1 k)^{-4}|p_j|}(p_j)=\emptyset,\forall\;i\neq j,
\end{equation}

\

\noindent or

\

\noindent (jj) there exist $0<r<R\leq 1$ such that
\begin{equation}\label{eq:jj1}
\int\limits_{A(r,R)} f \, dx\geq \frac{\tau_k}{(10 k)^2},
\end{equation}
and for any $N \in \N$, $N \geq 4(k+1)$, there exist $r\leq s_1<s_2<\cdots<s_{\scp{N+1}}\leq R$ such that
\begin{equation}\label{eq:jj2}
\int\limits_{A(s_i, s_{i+1})} f \, dx = \frac{1}{N} \int\limits_{A(r,R)}
f \, dx \; \qquad \quad \forall\,i=1,\ldots,N,
\end{equation}
and
\begin{equation}\label{eq:jj3}
\int\limits_{A(\frac{{s}}{C_1}, C_1 s)} f \, dx < \s_0, \qquad \quad \forall\,
s \in \left( C_1 r, \frac{R}{C_1} \right).
\end{equation}
\end{pro}

\

\noindent {\sc Proof.} We define for convenience
$$
\mathcal{A}_k=\{f\in L^{1}(B)\; | \; f > 0 \hbox{ a.e.}, \int_B f \, dx = 1,
 \, \, J_{k,\dt}(f) \leq 1-\tau_k\},
$$
and let
$$
\mathcal{A}_{k,0}=\left\{f\in L^{1}(B)\,|\, {(jj)} \; \mbox{holds}\;\mbox{for some}\;0<r<R \leq 1\right\}.
$$
For each $y\neq 0$ we denote $m_{\dt}(y;f)$ the integral
$$
m_{\dt}(y;f)=\int\limits_{B_{\d(10 C_1 k)^{-6}|y|}(y)} f\, dx.
$$
Consider the set
$$
\Lambda_k:=\left\{ \{x_1,\ldots,x_{k+1}\}\subset B \setminus \{0\}\, | \;
x_{i}\in B \setminus \bigcup\limits_{\ell\neq i}B_{\frac{\d |x_\ell|}{2}}(x_\ell),\;i=1,\ldots,k.\right\},
$$
and the number
$$
\s_k(\dt,\s_0):=\inf\limits_{f\in \mathcal{A}_k\setminus\mathcal{A}_{k,0}}
\sup\left\{ \min\limits_{i=1,\ldots,k+1}m_\dt(x_i;f)\,|\, \{x_1,\ldots,x_{k+1}\}\in \Lambda_k\right\}.
$$
A main step in our proof is the following

\

\noindent {\bf Claim:} $\s_k(\dt,\s_0)>0$.

\

\noindent {\sc Proof of the claim} Arguing by contradiction, for every $n\in\N$ there exists $f_n\in \mathcal{A}_k\setminus\mathcal{A}_{k,0}$ such that
\begin{equation}\label{H1-k}
\min\limits_{i=1,\ldots,k+1}m_\dt(x_i;f_n)\leq \frac{1}{n}\,,\;\;\forall \; \{x_1,\ldots,x_{k+1}\}\in \Lambda_k.
\end{equation}

\

\noindent For later use we fix here a positive number $0<\eps_0<<\frac{\tau_k}{(10 k)^2 10^k}$.
In the rest of this proof we will freely pass to subsequences which will not be relabelled and make use of the following:

\begin{lem}\label{lem1-k} Suppose that
\begin{equation}\label{Ha-1.1-k}
\int\limits_{A(r_{1,n},r_{2,n})} f_n\,dx\geq \frac{\tau_k}{(10 k)^2} >0,\; \qquad \quad
\forall n>\nu_0,
\end{equation}
for some $r_{1,n}<r_{2,n}$ and $\nu_0\in\N$. If there exists $0<\dt_0\leq \frac12$ such that
\begin{equation}\label{Ha-1.2-k}
\int\limits_{B_{\dt_0|x_n|}(x_n)}f_n\,dx\to 0,\;n\to+\i,\; \qquad \quad \forall x_n\in A(r_{1,n},r_{2,n}),
\end{equation}
then there exists $\nu_1\in\N$ such that $f_n\in \mathcal{A}_{k,0}$ for all $n>\nu_1$.
\end{lem}

\begin{pfn} {\sc of Lemma \ref{lem1-k}} We first prove that necessarily
\begin{equation}\label{L1}
\frac{r_{1,n}}{r_{2,n}}\to+\i,\;n\to+\i.
\end{equation}
We argue by contradiction and observe that then, up to the extraction of a subsequence, we could find $C>0$ such that
$$
A(r_{1,n},r_{2,n})\subseteq A(r_{1,n},C r_{1,n}).
$$
Observe that there exists $m=m(C)\in \N$ depending \un{only} on $C$ such that
$$
A(r,C r)\subseteq \bigcup\limits_{i=1}^{m} B_{\dt_0 |y_i|}(y_i),\;\{y_i\}_{i=1,\cdots,m}\subset A(r,C r).
$$
Therefore, by using \eqref{Ha-1.1-k}, \eqref{Ha-1.2-k} we obtain
$$
\frac{\tau_k}{(10 k)^2}\leq \int\limits_{A(r_{1,n},r_{2,n})} f_n\,dx\leq
\sum\limits_{i=1}^{m}\int\limits_{B_{\dt |x_{i,n}|}(x_{i,n})} f_n\,dx\leq m\,o(1),\;n\to+\i,
$$
which is the desired contradiction.

Next observe that $\int\limits_{A(s,t)} f_n\,dx$ is a continuous function of $s$ and $t$.
Hence, for any $N\in \N$ there exist $\{s_{1,n},s_{2,n},\cdots,s_{\scp{N+1},n}\}$ such that
$$
r_{1,n}\leq s_{1,n} \leq \cdots \leq s_{\scp{N+1},n}\leq r_{2,n},
$$
and
$$
\int\limits_{A(s_{i,n}, s_{i+1,n})}f_n\, dx=
\frac{1}{N} \int\limits_{A(r_{1,n},r_{2,n})} f_n\, dx\geq \frac{\frac{\tau_k}{(10 k)^2}}{N},\; \qquad \quad \forall\,i=1,\cdots,N.
$$
If for some $i\in\{1,\cdots,N+1\}$, along a subsequence we had
$$
\int\limits_{A(\frac{{s}_{i,n}}{C_1}, C_1 s_{i,n})}f_n\, dx\geq \s_0,\; \qquad \quad
\forall\, n\in\N
$$
then, by applying \eqref{L1} on $A(\frac{{s}_{i,n}}{C_1}, C_1 s_{i,n})$, we would obtain
$$
(C_1)^2=\frac{C_1 s_{i,n}}{\frac{{s}_{i,n}}{C_1}}\to+\i, \qquad \quad  n\to\i,
$$
which is the desired contradiction. Therefore there exists $\nu_1\in\N$ such that $f_n\in \mathcal{A}_{k,0}$
for any $n>\nu_1$. \end{pfn}

\

\noindent {\sc proof of the claim continued} There is no loss of generality in assuming
\begin{equation}\label{Ha-2-k.1}
m_\dt(x_{k+1,n},f_n)=\min\limits_{i=1,\ldots,k+1}m_\dt(x_i;f_n).
\end{equation}
Clearly \eqref{H1-k}, \eqref{Ha-2-k.1} and the definition of $\sigma_k(\d , \s_0)$
imply that for any $n\in\N$
\begin{equation}\label{Ha-3-k}
m_{\dt}(x_{k+1};f_n)\leq \frac{1}{n},\qquad \quad \forall x_{k+1}\neq 0 \hbox{ s.t. }
x_{k+1}\in B \setminus \bigcup\limits_{i=1,\ldots,k}B_{\frac{\d |x_{i,n}|}{2}}(x_{i,n}),
\end{equation}
with $\{x_{1,n},\cdots,x_{k,n},x_{k+1}\}\in\Lambda_k$. Set
$$
R_{n,(-)}=\min\limits_{i=1,\ldots,k}|x_{i,n}|\left(1-\frac{\dt}{2}\right),\qquad \quad
r_{n,(+)}=\max\limits_{i=1,\ldots,k}|x_{i,n}|\left(1+\frac{\dt}{2}\right),\;\forall\,n\in\N,
$$
and pick $0<r_{n,(-)}<R_{n,(-)}$, $r_{n,(+)}<R_{n,(+)}$ and $\nu\in\N$ such that
$$
\int\limits_{B_{r_{n,(-)}}} f_n\, dx<\frac{\eps_0}{2},
\;\int\limits_{B \setminus B_{R_{n,(+)}}} f_n\,dx<\frac{\eps_0}{2},\; \qquad \quad \forall\,n>\nu.
$$
If either $\int\limits_{A(r_{n,(-)},R_{n,(-)})} f_n\, dx\geq \frac{\tau_k}{(10 k)^2}$ or
$\int\limits_{A(r_{n,(+)},R_{n,(+)})} f_n\, dx\geq \frac{\tau_k}{(10 k)^2}$
for all $n>\nu_0$ for some $\nu_0\in\N$,
since \eqref{Ha-3-k} ensures that \eqref{Ha-1.2-k} holds on both $A(r_{n,(-)},R_{n,(-)})$ and
$A(r_{n,(+)},R_{n,(+)})$, then Lemma \ref{lem1-k} implies
that $f_n\in \mathcal{A}_{k,0}$ for all $n>\nu_1$, which is a contradiction.\\

Therefore, passing to a further subsequence if necessary, we can assume that
$$
\int\limits_{A(r_{n,(-)},R_{n,(-)})} f_n\, dx<\frac{\tau_k}{(10 k)^2}, \qquad
\hbox{ and } \qquad
\int\limits_{A(r_{n,(+)},R_{n,(+)})} f_n\, dx<\frac{\tau_k}{(10 k)^2},
$$
for any $n\in\N$, and in particular
$$
\int\limits_{A(R_{n,(-)},r_{n,(+)})} f_n\, dx\geq 1-2\frac{\tau_k}{(10 k)^2}-\eps_0>1-\frac{\tau_k}{10 k},
\qquad \quad   \forall\,n\in\N.
$$
Hence we conclude that
\begin{eqnarray}\label{last-k} \nonumber
  \int\limits_{A(R_{n,(-)},r_{n,(+)})\setminus \bigcup\limits_{i=1,\ldots,k}B_{\frac{\dt|x_{i,n}|}{2}}(x_{i,n})}
  f_n\, dx & \geq & \int\limits_{A(R_{n,(-)},r_{n,(+)})\setminus \bigcup\limits_{i=1,\ldots,k}B_{\dt|x_{i,n}|}(x_{i,n})} f_n\, dx
  \\ & \geq & 1-\frac{\tau_k}{10 k}-(1-\tau_k)>\frac{\tau_k}{2},
\end{eqnarray}
for all $n\in \N$. On the other hand, we have the following:

\

\begin{lem}\label{l:b3} There exists $\tilde{C} \geq 1$ such that
  $$
  \frac{r_{n,(+)}}{R_{n,(-)}}\leq \tilde{C},\qquad \quad \forall\,n\in\N.
  $$
\end{lem}

\begin{pfn} {\sc of Lemma \ref{l:b3}}
If the claim were false then, up to a relabelling of the indices,
we could find at least one index $i=i_n\in\{1,\ldots,k-1\}$ such that
$$
|x_{\ell,n}|\leq|x_{i_n,n}|<|x_{i_n+1,n}|\leq |x_{m,n}|,\; \qquad \quad
\forall \ell\leq i_n, \;\forall m\geq i_n+1,
$$
with the property that, passing to a subsequence if necessary,
\begin{equation}\label{04-12.1}
\lim\limits_{n\to+\i} \frac{|x_{i_n+1,n}|}{|x_{i_n,n}|}=+\i.
\end{equation}
Set
$$
R_{n,0}=|x_{i_n+1,n}|\left(1-\frac{\dt}{2}\right), \quad \;r_{n,0}=|x_{i_n,n}|\left(1+\frac{\dt}{2}\right), \qquad \quad \forall\,n\in\N.
$$
If $\int\limits_{A(r_{n,0},R_{n,0})} f_n\, dx\geq \frac{\tau_k}{(10 k)^2}$ for all $n>\nu_0$ for some $\nu_0\in\N$, since \eqref{Ha-3-k} and \eqref{04-12.1} together ensure that \eqref{Ha-1.2-k} holds on $A(r_{n,0},R_{n,0})$, then once more Lemma
\ref{lem1-k} implies that $f_n\in \mathcal{A}_{k,0}$ for all $n>\nu_1$, which is the desired contradiction.
\end{pfn}

 \

\noindent {\sc end of the proof of the claim} We are going to use Lemma \ref{l:b3}
together with \eqref{last-k} to obtain a contradiction. In fact, observe that there exists $\ell=\ell(\dt,\tilde{C})\in \N$ depending only on $\dt$ and $\tilde{C}$ such that
$$
A(R_{n,(-)},r_{n,(+)})\setminus \bigcup\limits_{i=1,\ldots,k}B_{\frac{\dt|x_{i,n}|}{2}}(x_{i,n})
\subseteq \bigcup\limits_{i=1}^{\ell}
B_{\dt(10C_1 k)^{-6}|y_{i,n}|}(y_{i,n}),
$$
where
$$
\{y_{i,n}\}_{i=1,\cdots,\ell}\subset
A(R_{n,(-)},r_{n,(+)})\setminus \bigcup\limits_{i=1,\ldots,k}B_{\frac{\dt|x_{i,n}|}{2}}(x_{i,n}).
$$
Hence \eqref{Ha-3-k} and \eqref{last-k} imply
$$
\frac{\tau_k}{2}\leq \int\limits_{A(R_{n,(-)},r_{n,(+)})\setminus
\bigcup\limits_{i=1,\ldots,k}B_{\frac{\dt|x_{i,n}|}{2}}(x_{i,n})} f_n\,dx\leq
\sum\limits_{i=1}^{\ell} m_{\dt}(y_{i,n};f_n)\leq \ell \,o(1),\qquad \quad
n\to+\i,
$$
which is the desired contradiction.\quad \rule{2mm}{2mm}

\

\noindent \noindent {\sc end of the proof of Proposition \ref{l:altprel}}
We claim that $\s=\min\{\frac{\s_k(\dt,\s_0)}{2},\s_0\}$ with $\dt=\dt_k=8(10 C_1 k)^{-4}$ satisfies the required  properties. In fact,
let us first assume that
$$
\frac{\s_k(\dt_k,\s_0)}{2}\leq \s_0,\quad\mbox{that is}\quad \s=\frac{\s_k(\dt_k,\s_0)}{2}.
$$
In this case, if $f$ does not satisfy $(jj)$, then by definition of $\s_k(\dt_k,\s_0)$
there exist $\{x_1,\cdots,x_{k+1}\}\in\Lambda_k$, such that
\begin{equation}\label{Dim-k}
\int\limits_{B_{(10 C_1 k)^{-8}|x_i|}(x_i)} f\, dx\geq
\int\limits_{B_{\d_k(10 C_1 k)^{-6}|x_i|}(x_i)} f\, dx\geq
\frac{\s_k(\dt_k,\s_0)}{2}=\s, \qquad \;i=1,\ldots,k+1.
\end{equation}
Next, let us prove that
\begin{equation}\label{Dim-2}
B_{\frac{\d |x_\ell|}{8}}(x_\ell)\cap B_{\frac{\d |x_m|}{8}}(x_m)=\emptyset,\qquad
\quad \forall \{\ell,m\}\subset\{1,\ldots,k\}, \quad \ell\neq m.
\end{equation}
If $|x_\ell|\leq 2|x_m|$ and $x\in B_{\frac{\dt|x_{m}|}{8}}(x_{m})$, then
$$
d_g(x,x_\ell)>\frac{\dt|x_m|}{2}-\frac{\dt|x_m|}{8}=\frac{3\dt|x_m|}{8}>\frac{\dt|x_m|}{4}\geq \frac{\dt|x_\ell|}{8},
$$
while if $|x_\ell|> 2|x_m|$ and $x\in B_{\frac{\dt|x_{m}|}{8}}(x_{m})$, then
$$
d_g(x,x_\ell)>|x_\ell|-|x_m|-\frac{\dt|x_m|}{8}>\frac{\dt|x_m|}{8},
$$
that is,  \eqref{Dim-2} holds.
Hence, if $\dt=\dt_k$ we see that $(j)$ is satisfied with $\{p_1,\cdots,p_{k+1}\}=\{x_1,\cdots,x_{k+1}\}$ and the desired property holds with
$\s\leq \s_0$.

\

\noindent On the other hand, if
$$
\frac{\s_k(\dt_k,\s_0)}{2}> \s_0,\quad\mbox{that is}\quad \s=\s_0,
$$
and $(jj)$ is not satisfied then by definition of $\s_k(\dt_k,\s_0)$ we can find
$\{x_1,\cdots,x_{k+1}\}$ as above such that \eqref{Dim-k} and \eqref{Dim-2} with $\dt=\dt_k$ hold,
so that $(j)$ is satisfied with $\{p_1,\cdots,p_{k+1}\}=\{x_1,\cdots,x_{k+1}\}$ and $\s=\s_0$.
This concludes the proof. \quad \rule{2mm}{2mm}

\

\noindent In the next subsections we prove Proposition \ref{p:impr} in both alternatives
of Proposition \ref{l:altprel} choosing $\sigma_0$ as
\begin{equation}\label{eq:epsilon}
    \s_0 = \frac{\t_k}{100 k^2} \frac{\e}{4 (k+1) \log C_1}.
\end{equation}

\

\subsection{Proof of Proposition \ref{p:impr} in case $(jj)$}\label{ss:van}

We will argue that if a certain fixed amount of conformal volume is
{\em diluted} in a large portion of the cylinder, then we  can divide it
into $N$ parts, with $N$ large enough, and choose the one with smallest
Dirichlet energy for using Lemma \ref{l:ext}.

Letting $r, R$ be as in $(jj)$,  we choose a large number $N$, depending
on $k_\a$ and $\e$, such that
$$
  N = \left[ \frac{8(k_\a+1)}{\e} \right],
$$
where the square bracket stands for the integer part.

We next choose $r \leq s_1 < \dots < s_{N+1} \leq R$
such that \eqref{eq:jj2} and \eqref{eq:jj3} hold. We notice immediately that
by the choice of $N$ one has
\begin{equation}\label{eq:inttildesi}
      \int_{A(s_i, s_{i+1})} \tilde{f}_u dx \geq
      \frac{\e}{16(k_\a+1)} \frac{\t_k}{(10k)^2}; \qquad \quad i = 1, \dots, N.
\end{equation}
We also claim that for every index $i$ the intersection $A\left(\frac{s_i}{C_1}, C_1 s_i
\right) \cap A \left(\frac{s_{i+1}}{C_1}, C_1 s_{i+1}\right)$ is empty. In fact,
if this were not the case we would have by \eqref{eq:jj3} and \eqref{eq:epsilon}
\begin{eqnarray*}
  \frac{\e}{16(k_\a+1)} \frac{\t_k}{(10k)^2} & \leq & \int_{A(s_i, s_{i+1})} \tilde{f}_u
     dx \leq \int_{A\left(\frac{s_i}{C_1}, C_1 s_i
     \right)} \tilde{f}_u dx + \int_{A\left(\frac{s_{i+1}}{C_1}, C_1 s_{i+1} \right)}
     \tilde{f}_u dx \\ & \leq & \frac{\e}{16(k_\a+1)} \frac{2 \t_k}{(10k)^2 \log C_1},
\end{eqnarray*}
which is a contradiction by the choice of $C_1$.
%
%
%
We can now choose an index $i$ for which
\begin{equation}\label{eq:indinu}
    \int_{A(s_i/C_1, C_1 s_{i+1})}
    |\n u|^2 dx \leq \frac{2}{N} \int_B |\n u|^2 dx \leq
    \frac{\e}{4(k_\a+1)} \int_B |\n u|^2 dx.
\end{equation}
We can also choose $\tilde{s}_i \in \left[ \frac{2s_i}{C_1}, \frac{s_i}{2} \right]$
and $\tilde{s}_{i+1} \in \left[ 2 s_{i+1}, \frac{1}{2} C_1 s_{i+1} \right]$
such that
$$
  \int_{A(\tilde{s}_i/2, 2 \tilde{s}_i)} |\n u|^2 dx \leq \frac{4}{\log C_1}
  \int_{A(s_i/C_1, C_1 s_{i+1})} |\n u|^2 dx;
$$
$$
  \int_{A(\tilde{s}_{i+1}/2, 2 \tilde{s}_{i+1})} |\n u|^2 dx \leq \frac{4}{\log C_1}
    \int_{A(s_i/C_1, C_1 s_{i+1})} |\n u|^2 dx,
$$
so by \eqref{eq:indinu} we have
  \begin{equation}\label{eq:kkk1}
  \int_{A(\tilde{s}_i/2, 2 \tilde{s}_i)} |\n u|^2 dx \leq \frac{4}{\log C_1}
   \frac{\e}{4(k_\a+1)} \int_B |\n u|^2 dx;
  \end{equation}
  \begin{equation}\label{eq:kkk2}
  \int_{A(\tilde{s}_{i+1}/2, 2 \tilde{s}_{i+1})} |\n u|^2 dx \leq \frac{4}{\log C_1}
    \frac{\e}{4(k_\a+1)} \int_B |\n u|^2 dx.
  \end{equation}
We define then a new function $\tilde{u}$ as the harmonic lifting of
$u$ inside $B_{\tilde{s}_i}(0)$: recalling the definition in \eqref{eq:hlift},
we set
$$
  \tilde{u}(x) = \mathcal{H}_{B_{\tilde{s}_i}}(u).
$$
By Lemma \ref{l:ext}, \eqref{eq:indinu} and \eqref{eq:kkk1} we have that
\begin{equation}\label{eq:dinternomod}
    \int_{B_{2 \tilde{s}_{i+1}}} |\n \tilde{u}|^2 dx \leq \frac{\e}{2(k_\a+1)}
    \left( 1 + \frac{4 C_0}{\log C_1} \right) \int_B |\n u|^2 dx.
\end{equation}
We then apply Proposition \ref{l:mr} with $s = \tilde{s}_{i+1}/2$,
$\eta = \frac{4}{\log C_1} \frac{\e}{4(k_\a+1)}$ and $\t =
\frac{\e}{16(k_\a+1)} \frac{\t_k}{(10k)^2}$ (see \eqref{eq:inttildesi}, \eqref{eq:kkk2}
and the choice of $\tilde{s}_{i+1}$)  and use \eqref{eq:dinternomod} to find
\begin{eqnarray*}
  (1+\a)\log \int_B f_u dx & \leq & \frac{1}{4 \pi} \left( \a \int_{B_{\tilde{s}_{i+1}/2}}
  |\n u|^2 dx + \int_{B \setminus B_{\tilde{s}_{i+1}/2}} |\n u|^2 dx
  + C_0 \eta \int_B |\n u|^2 dx \right) + \tilde{C} \\
  & \leq & \frac{1}{4 \pi} \left(  1 +  \frac{\a \e}{2(k_\a+1)} \left(
  1 + \frac{4 C_0}{\log C_1} \right) + C_0 \frac{4 \e}{4(k_\a+1) \log C_1} \right)
  \int_B |\n u|^2 dx + \tilde{C} \\
  & \leq & \left( 1 + \frac{\a \e}{2(k_\a+1)} +
  C_0 \frac{4 \e (\a + 1)}{4(k_\a+1) \log C_1} \right)  \int_B |\n u|^2 dx + \tilde{C},
\end{eqnarray*}
where $\tilde{C}$ depends on $\a$ and  $\e$. By the choice of $C_1$ the last formula
implies
$$
   (1+\a)\log \int_B f_u dx \leq \frac{1}{4 \pi} \left( 1 + \e \right)
  \int_B |\n u|^2 dx + \tilde{C},
$$
and in turn
$$
   \min \{ 1 + k, 1 + \a \} \, \log \int_B f_u dx \leq \frac{1}{4 \pi}
   \left( 1 + \e \right) \int_B |\n u|^2 dx + \tilde{C},
$$
which concludes the proof.

\

\subsection{Proof of Proposition \ref{p:impr} in case $(j)$}\label{ss:impro}

If $p_1, \dots, p_{k+1}$ are as in \eqref{eq:j1}, then
there exist $\theta_i$, $i = 1, \dots, k+1$ such that
\begin{equation}\label{eq:thi}
    (10 C_1 k)^{-6} |p_i| \leq \theta_i \leq (10 C_1 k)^{-5} |p_i|; \qquad \quad
    \int_{B_{4 \th_i}(p_i) \setminus B_{\th _i}(p_i)} |\n u|^2 dx
    \leq \eta \int_B |\n u|^2 dx,
\end{equation}
where
$$
  \eta = \frac{1}{\log (10 k C_1)}.
$$
We can also assume that either $B_{4 \th_i}(p_i) \subseteq B$ or that
$d(p_i, \pa B) \leq \frac{1}{4} \th_i$ (we require these conditions in view of
Remark \ref{r:centro}).

We next select an index $\ov{i}$ such that
$$
  \ov{D}:= \frac{1}{4 \pi} \int_{B_{\th_{\ov{i}}(p_{\ov{i}})}} |\n u|^2 dx =
  \min_{i \in \{1, \dots, k+1 \}} \frac{1}{4 \pi} \int_{B_{\th_i}(p_{i})} |\n u|^2 dx,
$$
and then another index $\tilde{i}$ for which
$$
  \tilde{D}:= \frac{1}{4 \pi} \int_{B_{\th_{\tilde{i}}(p_{\tilde{i}})}} |\n u|^2 dx =
    \min_{i \neq \ov{i}} \frac{1}{4 \pi} \int_{B_{\th_i}(p_{i})} |\n u|^2 dx.
$$
Below, we  set for convenience
$$
  D_1 = \frac{1}{4 \pi} \int_{\cup_{i=1}^{k+1} B_{\th_i}(p_{i})} |\n u|^2 dx;
  \qquad \quad D_2 = \frac{1}{4 \pi} \int_{B \setminus \cup_{i=1}^{k+1}
  B_{\th _i}(p_i)} |\n u|^2 dx,
$$
and
$$
   D = \frac{1}{4 \pi} \int_B |\n u|^2 dx = D_1 + D_2.
$$
Notice that, by our choices of $\ov{i}$ and $\tilde{i}$, $\ov{D}$
and $\tilde{D}$ satisfy
\begin{equation}\label{eq:ovtd}
  \ov{D} \leq \frac{D_1}{k+1}; \qquad \quad \qquad \tilde{D} \leq \frac{D_1 - \ov{D}}{k}.
\end{equation}
We then consider a modified function $\check{u}$ defined as
$$
  \check{u} = \mathcal{H}_{\cup_{i\neq \ov{i},\tilde{i}} B_{\th_{i}}(p_{i})}(u).
$$
Notice that, by construction and by Lemma \ref{l:ext}, one has
\begin{equation}\label{eq:estentildeuu}
  \frac{1}{4 \pi} \int_{B \setminus B_{\th_{\ov{i}}(p_{\ov{i}})}} |\n \check{u}|^2 dx
  \leq \tilde{D} + D_2 +  C_0 k \eta D, 
\end{equation}
and also, by \eqref{eq:j1}
\begin{equation}\label{eq:estinttildeuu}
  \log \int_{B} f_u dx \leq \log \int_{B_{\th_{\ov{i}}}} f_{\check{u}} 
  dx + C_\s; \qquad \quad  \log \int_{B} f_u dx \leq \log 
  \int_{B \setminus B_{4 \th_{\ov{i}}}} f_{\check{u}} dx + C_\s. 
\end{equation}
Using then \eqref{eq:mmmm} and \eqref{eq:mmmm2222} for $\check{u}$, with $B_s(p_{\ov{i}})$,
$s = \th_i$, $\t = \s$, and taking Remark \ref{r:centro} into account (in which we 
allow $\ov{C}$ to depend on $\e$, $\al$ but not on the $p_i$'s), from 
\eqref{eq:estinttildeuu} we get 
\begin{equation}\label{eq:mmmm-mod}
   \log \int_{B} f_{u} dx \leq 
  \frac{1}{4 \pi} \left( \int_{B}
   |\n \check{u}|^2 dx + 2 C_0 \eta \int_B |\n u|^2 dx \right) + \tilde{C}, 
\end{equation}
where $\tilde{C}$ depends on $\eta, \t$ and $\a$. From \eqref{eq:estentildeuu} 
we obtain 
$$
   (\a + 1) \log \int_B f_u dx \leq \a \ov{D} + \tilde{D} + D_2 +
   C_0 (k+2) \eta D + \tilde{C}. 
$$
Then by the second inequality in \eqref{eq:ovtd} one finds
$$
  (\a + 1) \log \int_B f_u dx \leq \a \ov{D} + \frac{D_1 - \ov{D}}{k}
  + D_2 +  C_0 (k+2) \eta D + \tilde{C},
$$
which implies
\begin{equation}\label{eq:finfin}
   (\a + 1) \log \int_B f_u dx \leq \frac{1}{k} D_1 + \left( \a - \frac 1k \right)
       \ov{D} + D_2 +  C_0 k \eta D + \tilde{C}.
\end{equation}
If $\a \leq 1$ then necessarily $k = 1$, so the coefficient of $\ov{D}$ in the
latter formula is negative and can be discarded, yielding
$$
 (\a + 1) \log \int_B f_u dx \leq D_1 + D_2 +  C_0 (k+2) \eta D + \tilde{C}
 \leq  D +  C_0 k \eta D + \tilde{C} \leq (1+\e) D + \tilde{C}
$$
which gives the conclusion, by our choices of $\eta$ and $C_1$. 

On the other hand, if $\a > 1$ we have that $\a - \frac 1k > 0$ and hence,
since $\ov{D} < \frac{1}{k+1} D_1$ (see the first inequality in \eqref{eq:ovtd}),
\eqref{eq:finfin} gives
\begin{eqnarray*}
  (\a + 1) \log \int_B f_u dx & \leq & \left( \frac 1k + \frac{\a k - 1}{k}
    \frac{1}{k+1} \right) D_1 + D_2 +  C_0 k \eta D + \tilde{C} \\
    & \leq & \frac{k+1 + \a k -1}{k(k+1)}  D_1 + D_2 +  C_0 k \eta D + \tilde{C} \\
     & \leq & \frac{\a + 1}{k+1}  D_1 + D_2 +  \e D + \tilde{C}.
\end{eqnarray*}
If $\a \leq k$ this implies
$$
  (\a + 1) \log \int_B f_u dx \leq  D_1 + D_2 +  \e D \leq (1 + \e) D + \tilde{C},
$$
as desired.

If instead $k < \a$ we obtain
$$
  (k + 1) \log \int_B f_u dx \leq  D_1 + \frac{k + 1}{\a+1} D_2 +  \e D
  \leq (1+\e) D + \tilde{C},
$$
which still gives the conclusion.

\

\

\noindent From the latter proposition we immediately deduce
the following lower bound on $I_\rho$, which can be obtained
choosing $\e > 0$ small enough.

\begin{cor}\label{c:lowbd}  Let $\d_k$ and $\t_k$ be so small that Proposition
\ref{p:proj} applies, and let $\d \leq \d_k$. Let $k \in \{1, \dots, k_\a\}$: then
there exists a constant $C_{k,\a}$, depending only on $k$ and $\a$, such that
$$
  I_{\rho, \underline{\a}}(u) \geq - C_{k,\a}
$$
for all functions $u$ such that $J_{k,\d}(\tilde{f}_u) \leq 1 - \t_k$.
\end{cor}

\noindent As a consequence of the last corollary and of Proposition
\ref{p:homeo} we obtain an explicit
condition which guarantees lower bounds on $I_{\rho, \underline{\a}}$.

\begin{cor}\label{c:moments} Let $\d_k$ and $\t_k$ be so small that Proposition
\ref{p:proj} applies, and let $\d \leq \d_k$. Let $k \in \{1, \dots, k_\a\}$, and
let $F_k$ denote the map in \eqref{eq:FkFk}. Then there exists a constant
$C_{k,\a}$, depending only on $k$ and $\a$, such that
$$
  I_{\rho, \underline{\a}}(u) \geq - C_{k,\a}
$$
provided $F_k(\tilde{f}_u) = 0$.
\end{cor}

\section{Proof of the existence  and non existence results}\label{s:ex}

\noindent In this section we provide applications of the improved inequality in
Proposition \ref{p:impr} to the existence of solutions to \eqref{eq:e-1}. We give
full details in two simple cases, namely in the unit
ball with Dirichlet boundary data and one singularity, as well as on the sphere with two
singularities, see Remark \ref{r:ab} for more general situations. The variational argument combines different known strategies, therefore we will be quite sketchy in some parts.

We then prove one may have non existence of solutions in case the
assumptions on $\rho$ are dropped, showing that the hypotheses of
Theorems \ref{t:exdisk} and \ref{t:sphere} are sharp.

\

\subsection{Proof of Theorem \ref{t:exdisk} (for $m = 1$ in simply connected domains)}

First of all, through a Riemann map we can reduce ourselves to
the case of the unit ball $B$ with the singularity at the origin.
We let $k$ be the unique integer for which $\rho \in (4 k \pi,
4 (k+1)\pi)$, and we let $F_k$ denote the map  in \eqref{eq:FkFk},
which realizes a homeomorphism between $(S^1)_k$ and $\mathcal{S}_k$,
see Proposition \ref{p:homeo}.

Choose a non negative cut-off function $\chi$ such that
$$
  \begin{cases}
  \chi \in C^\infty_c(B); & \\
  \chi(x) \equiv 1 & \hbox{ in } B_{\frac{3}{4}},
  \end{cases}
$$
and for $\s = \sum_{i=1}^k t_i \d_{\th_i} \in (S^1)_k$, $\l > 0$, we define
the test function
\begin{equation}\label{eq:testdisk}
   \var_{\l ,\s}(x) = \chi(x) \log \sum_{i=1}^{k} t_i \left(
   \frac{\l}{1 + \l^2 \left| y - \frac{1}{2} x_i \right|^2} \right)^2,
   \qquad \quad x_i = (\cos \th_i, \sin \th_i).
\end{equation}
Reasoning as in \cite{dj} (see also \cite{mald} for a simpler proof of this
estimate) one can obtain the following result with minor modifications of the proof.

\begin{lem}\label{l:lowen}
Let  $\var_{\l ,\s}$ be defined as in \eqref{eq:testdisk}. Then as
$\l \to + \infty$ one has
   \begin{equation}\label{eq:flatconv}
   d_{K-R}(\tilde{f}_{\var_{\l ,\s}}, \tilde{\s}) \to 0, \qquad \quad
   \tilde{\s} = \sum_{i=1}^k t_i \d_{\frac{1}{2} x_i},
   \end{equation}
and
$$
  I_{\rho, \underline{\al}}(\var_{\l ,\s}) \to - \infty
$$
uniformly for $\s \in (S^1)_k$.
\end{lem}

\

\noindent We next define the variational scheme which will allow us to
find existence of solutions. Recalling that $\mathcal{U}_k$ denotes the interior
of $\mathcal{S}_k$ in $\C^k$,  consider the family of continuous maps
$$
  \mathcal{K}_{\l,\rho} = \left\{ \mathfrak{h} : \mathcal{U}_k \to H^1_0(B)
  \; : \; \mathfrak{h}(y) = \var_{\l,F_k^{-1}(y)} \hbox{ for every } y \in
  \mathcal{S}_k = \pa \mathcal{U}_k \right\}.
$$
We define also the min-max value
$$
  \ov{\mathcal{K}}_{\l,\rho} = \inf_{\mathfrak{h} \in \mathcal{K}_{\l,\rho}}
  \sup_{z \in \mathcal{U}_k} I_{\rho, \underline{\a}}(\mathfrak{h}(z)).
$$
We have then the following result, which implies the conclusion of
Theorem \ref{t:exdisk}.

\begin{pro}\label{p:crit0} Under the assumptions of Theorem \ref{t:exdisk},
if $\l$ is sufficiently large then
$$
   \ov{\mathcal{K}}_{\l,\rho} > \sup_{y \in \mathcal{S}_k}
   I_{\rho, \underline{\a}}(\var_{\l,F_k^{-1}(y)}).
$$
Moreover $\ov{\mathcal{K}}_{\l,\rho}$ is a critical value of $I_{\rho,\underline{\a}}$.
\end{pro}

\begin{pf} If $C := C_{k,\a}$ is as in Corollary \ref{c:moments}, we let $L = 4 C$, and
choose $\l$ to be so large that
$$
  \sup_{y \in \mathcal{S}_k} I_{\rho,\underline{\a}}(\var_{\l,F_k^{-1}(y)}) < - L,
$$
which is possible in view of Lemma \ref{l:lowen}.

We are going to show that $\ov{\mathcal{K}}_{\l,\rho} > - \frac{L}{2}$.
Indeed, assume by contradiction that there exists a continuous
$\mathfrak{h}_0$ such that
\begin{equation}\label{eq:h0}
 \mathfrak{h}_0 \in \mathcal{K}_{\l,\rho} \qquad \hbox{ and } \qquad
 \sup_{z \in \mathcal{U}_k} I_{\rho,\underline{\a}}(\mathfrak{h}_0(z)) \leq - \frac 12 L.
\end{equation}
Then, by our choice of $L$, Corollary \ref{c:lowbd} and Proposition \ref{p:proj} would apply,
yielding a continuous map
$F_{\l,\rho} : \mathcal{U}_k \to \mathcal{S}_k$ defined as the
composition
$$
  F_{\l,\rho} = \Xi_k \circ \mathfrak{h}_0.
$$
Notice that, since $\mathfrak{h}_0 \in \mathcal{K}_{\l,\rho}$, $\mathfrak{h}_0(\cdot)$
coincides with $\var_{\l,F_k^{-1}(\cdot)}$ on $\mathcal{S}_k = \pa \mathcal{U}_k$,
so by \eqref{eq:flatconv} we deduce that
\begin{equation}\label{eq:hom}
      F_{\l,\rho}|_{\mathcal{S}_k}  \hbox{ is homotopic to } Id|_{\mathcal{S}_k}:
\end{equation}
the homotopy is obtained  by letting the parameter $\l$ tend to
$+\infty$. Since $\mathcal{S}_k$ is homeomorphic to $S^{2k-1}$, it is non
contractible, and we obtain a contradiction to \eqref{eq:hom}. This proves
$ \ov{\mathcal{K}}_{\l,\rho} > \sup_{y \in \mathcal{S}_k}
   I_{\rho,\underline{\a}}(\var_{\l,F_k^{-1}(y)})$.

\

\noindent To check that $\ov{\mathcal{K}}_{\l,\rho}$ is a critical
level is rather standard, as one can use a monotonicity method
from \cite{lucia}, \cite{str}. Consider a sequence $\rho_n \to \rho$
and the corresponding functionals $I_{\rho_n,\underline{\a}}$. All the above estimates,
including also those from the previous sections, can be worked out for $I_{\rho_n,\underline{\a}}$
as well with minor changes, if $n$ is large enough.

We then define the min-max value $\tilde{\mathcal{K}}_{\l,\rho} :=
\frac{\ov{\mathcal{K}}_{\l,\rho}}{\rho}$, which corresponds to the
functional $\frac{I_{\rho,\underline{\a}}}{\rho}$. It is immediate to see that
$$
  \rho \mapsto \tilde{\mathcal{K}}_{\l,\rho} \qquad \hbox{ is monotone},
$$
and that, reasoning as in \cite{djlw}, there exists a subsequence of
$(\rho_n)_n$ such that $I_{\rho_n,\underline{\a}}$ has a solution $u_n$ at level
$\ov{\mathcal{K}}_{\l,\rho_n}$. Then, applying Corollary \ref{c:comp}
and passing to a further subsequence, we obtain that $u_n$
converges to a critical point $u$ of $I_{\rho,\underline{\a}}$ at level
$\ov{\mathcal{K}}_{\l,\rho}$.
\end{pf}

\subsection{Proof of Theorem \ref{t:sphere} (for $m = 2$)}

\noindent The argument is very similar in spirit to the
previous case. We list the main changes which are necessary
to deal with this situation, especially for what concerns the improved
Moser-Trudinger inequality.

\

\noindent First of all, using again a M\"obius map on $S^2$, we can assume that
the two singularities $p_1$ and $p_2$ are antipodal, and coincide
respectively with the south and the north pole of $S^2$, viewed as
the standard sphere embedded in $\R^3$.

Given a small $\d > 0$ we can define the following quantity, analogous
to the one in \eqref{eq:intk}
\begin{equation}\label{eq:tintk}
    \tilde{J}_{k,\d}(\tilde{f}_u) = \sup_{x_1, \dots, x_k \neq \{p_1, p_2\}} \int_{\cup_{i=1}^k
    B_{\d \min \{ d(x_i, p_1), d(x_i, p_2) \}(x_i)}} \tilde{f}_u dV_g,
\end{equation}
as well as the measure on the unit circle (viewed as the $(x,y)$ plane in
$\R^3$ intersected with $S^2$)
$$
  \tilde{\mu}_f(A) = \int_{\tilde{\pi}^{-1}(A)} \tilde{f}_u dV_{g};  \qquad
  \quad A \subseteq S^1,
$$
where $\tilde{\pi} : S^2 \setminus \{ p_1, p_2 \} \to S^1$ stands for the projection
onto the equator along the meridians.

Reasoning as in Section 3, if $\tilde{J}_{k,\d}(\tilde{f}_u) > 1 - \t_k$, $\d \leq \d_k$,
we can project continuously $u$ onto the $k$-barycenters of $S^1$, $(S^1)_k$.
For the case $\tilde{J}_{k,\d}(\tilde{f}_u) \leq 1 - \t_k$ we have a counterpart of
Proposition \ref{p:impr}.

\begin{pro}\label{p:imprsph} Let $\e > 0$, and let $k \in \{1, \dots, k_\a\}$.
Let $\d_k$ and $\t_k$ be so small that Proposition \ref{p:proj}
applies, and let $\d \leq \d_k$. Then there exists a constant
$C_{\e,\a_1, \al_2}$, depending only on $\e$, $\a_1$ and $\a_2$,
such that
$$
  \log \int_{S^2} f_u dV_{g} \leq \frac{1+\e}{4\pi \min\{ 1 + k, 1+\a_1, 1 + \al_2 \}}
  \int_{S^2} |\n u|^2 dV_g + C_{\e,\a_1, \al_2} + 2 \fint_{S^2} u \, dV_{g}
$$
for all functions $u$ verifying
$$
  \tilde{J}_{k,\d}(\tilde{f}_u) \leq 1 - \t_k.
$$
\end{pro}

\noindent To check this statement, one can reason as in the proof
of Proposition \ref{p:impr}, with two main differences. The first
is that the average of $u$ on $S^2$ should be added to the
right-hand side of the inequality, since there are no boundary
data in this case (compare \eqref{eq:mtdisk} and
\eqref{eq:mtdisk0}): it will not be a loss of generality to assume
that $\fint_{S^2} u \, dV_{g} = 0$. The second is that in case
$(j)$ (resp. in case $(jj)$) the points $x_{i}$ (resp. the regione
of {\em vanishing} for the measure $\tilde{f}_u$) can lie near
either $p_1$ or $p_2$.

Suppose that $(jj)$ holds, and let $p_i$ be a point near which
vanishing occurs. Then the previous arguments yield the inequality 
$$
\min\{ 1 + k, 1+\a_i\} \log \int_{S^2} f_u dV_g \leq
\frac{1+\e}{4\pi} \int_{S^2} |\n u|^2 dV_g + C_{\e,\a_1,\a_2},
$$
which implies
$$
  \min\{ 1 + k, 1+\a_1, 1 + \al_2 \} \log \int_{S^2} f_u dV_g \leq
  \frac{1+\e}{4\pi} \int_{S^2} |\n u|^2 dV_g
+ C_{\e, \a_1, \a_2},
$$
namely the desired conclusion.

If $(j)$ holds instead, there will be $k_1$ points among the
$x_i$'s approaching $p_1$, and $k_2$ points which either lie in a
fixed compact set of $S^2 \setminus \{p_1, p_2\}$ or approaching
$p_2$, with $k_1 + k_2 = 1 + k$. Applying Proposition \ref{p:impr}
(twice, with $B$ replaced by two spherical caps whose union is
$S^2$ and whose boundaries are well separated from the points
$x_i$) and Lemma \ref{l:medie} one finds that
$$
  \left( \min \{k_1, 1 + \al_1\} + \min \{k_2, 1 + \al_2\} \right)
  \log \int_{S^2} f_u dV_g \leq \frac{1+\e}{4\pi} \int_{S^2}
  |\n u|^2 dV_g + C_{\e,\a_1,\a_2}.
$$
Then it is enough to use the elementary inequality
$$
   \min\{ 1 + k, 1+\a_1, 1 + \al_2 \} \leq \min \{k_1, 1 + \al_1\}
   + \min \{k_2, 1 + \al_2\},
$$
to obtain again the conclusion.

\

\noindent For $\th_1, \dots, \th_k \in S^1$ and for
$\hat{\s} = \sum_{i=1}^k t_i \d_{\th_i}$, the following test function replaces
$\var_{\l ,\s}$ in \eqref{eq:testdisk}
\begin{equation}\label{eq:testdisk2}
   \hat{\var}_{\l ,\hat{\s}}(x) = \log \sum_{i=1}^{k} t_i \left(
   \frac{\l}{1 + \l^2 d(y, x_i)^2} \right)^2, \qquad \quad
   x_i = (\cos \th_i, \sin \th_i, 0) \in S^2.
\end{equation}
One can then prove the counterpart of Lemma \ref{l:lowen}.

\begin{lem}\label{l:lowen2}
Let  $\hat{\var}_{\l ,\hat{\s}}$ be defined as in \eqref{eq:testdisk2}. Then as
$\l \to + \infty$ one has
   \begin{equation}\label{eq:flatconv2}
   d_{K-R}(\tilde{f}_{\hat{\var}_{\l ,\hat{\s}}}, \hat{\s}) \to 0, \qquad \quad
  I_{\rho,\underline{\a}}(\hat{\var}_{\l ,\hat{\s}}) \to - \infty
\end{equation}
uniformly for $\hat{\s} \in (S^1)_k$.
\end{lem}

\

\noindent In the above lemma, with an abuse of notation, we are identifying $S^1$ as
the equator of $S^2$. Considering now the class of continuous maps
$$
  \hat{\mathcal{K}}_{\l,\rho} = \left\{ \mathfrak{h} : \mathcal{U}_k \to H^1_0(B)
  \; : \; \mathfrak{h}(y) = \hat{\var}_{\l,F_k^{-1}(y)} \hbox{ for every } y \in
  \mathcal{S}_k = \pa \mathcal{U}_k \right\}.
$$
and the min-max value
$$
  \hat{\ov{\mathcal{K}}}_{\l,\rho} = \inf_{\mathfrak{h} \in \hat{\mathcal{K}}_{\l,\rho}}
  \sup_{z \in \mathcal{U}_k} I_{\rho,\underline{\a}}(\mathfrak{h}(z)),
$$
one has the counterpart of Proposition \ref{p:crit0}.

\begin{pro}\label{p:crit02} Under the assumptions of Theorem \ref{t:exdisk},
if $\l$ is sufficiently large then
$$
   \hat{\ov{\mathcal{K}}}_{\l,\rho} > \sup_{y \in \mathcal{S}_k}
   I_{\rho,\underline{\a}}(\hat{\var}_{\l,F_k^{-1}(y)}).
$$
Moreover $\hat{\ov{\mathcal{K}}}_{\l,\rho}$ is a critical value of $I_{\rho,\underline{\a}}$.
\end{pro}

\

\begin{rem}\label{r:ab}  As anticipated in Remark \ref{r:commenti}, the above min-max
method can also be applied to the case of more singularities, multiply connected domains or to surfaces with positive genus, combining the present approach to the one in \cite{bdm}.

Regarding Theorem \ref{t:exdisk} for $m \geq 2$ or domain $\O$ of $\R^2$
which is not simply connected we argue as follows. Choosing an index $\ov{i}$
for which $\a_{\ov{i}} = \min_{i=1, \dots, m} \a_{i}$, one can find a simple
curve $\g$ in $\O \setminus \cup_i p_i$ non contractible in $\O \setminus
p_{\ov{i}}$ such that there is a continuous map $\G : \O \setminus
p_{\ov{i}} \to \g$ satisfying $\G|_{\g} = Id|_{\g}$.

On one hand, it is possible to associate to each $\tilde{f}_u$, $u \in H^1_0(\O)$,
a unit measure on $S^1$ via the push-forward of $\G$, and hence introduce a
counterpart of the map $F_k$. On the other, one can use for the min-max scheme a
test function as in \eqref{eq:testdisk}, but with the points $x_i$ distributed
on $\g$.

The combination of these two facts allows to repeat the above procedure: in the case
of the sphere with $m > 2$ or for compact surfaces with positive genus one can argue similarly.
\end{rem}

\

\subsection{A non existence result on the unit ball}

Here we show that our theorem on simply connected domains is sharp. Actually
we provide a sketchy proof of the following well known fact:\\

\begin{pro}\label{p:nonexdisk}
 If the following problem admits a solution $u\in H^{1}_0(B)$
 \begin{equation}\label{eq:Diri}
 \graf{
 -\Delta u=\rho \frac{|x|^{2\al} e^{2u}}{\int_{B}
         |x|^{2\al}e^{2u}dx} & \mbox{in}\quad B\\
 u=0 & \mbox{on}\quad \pa B,
 }
 \end{equation}
 then necessarily $\rho<4\pi(1+\al)$.
\end{pro}

\begin{pf}
Set $V(x)=\rho|x|^{2\al} \left(\int_{B}|x|^{2\al}e^{2u}dx\right)^{-1}$.
   By using the fact that $u=0$ on $\pa B$, then the Poho\u{z}aev
   identity for the equation in \eqref{eq:Diri} reads
   $$
   -\frac{1}{2}\int\limits_{\pa B}(x,\nu)(u_\nu)^2 d\s = \frac{1}{2}
   \int\limits_{\pa B}(x,\nu)V e^{2u}  d \s -\frac{1}{2}
   \int\limits_{B} \left[2V
   e^{2u}+ \langle x,\nabla \log{V} \rangle Ve^{2u}\right] dx,
   $$
   where $\nu$ is the unit outer normal to $\pa B$ and $u_\nu=(\nu,\nabla u)$. Since $(x,\nu)=1$ on $\pa B$,
   the Cauchy-Schwarz inequality then yields
\begin{eqnarray*}
    -\frac{1}{4\pi}\left(\,\,\int\limits_{\pa B}u_\nu \, d \s \right)^2 & = &
    -\frac{1}{2}\left(\,\,\int\limits_{\pa B}\frac{d \s}{(x,\nu)}\right)^{-1}
       \left(\,\,\int\limits_{\pa B}u_\nu \, d \s \right)^2 \geq
       -\frac{1}{2}\int\limits_{\pa B}(x,\nu)(u_\nu)^2 d\s \\ & = &
       \frac{1}{2}\int\limits_{\pa B}(x,\nu)V e^{2u} d \s -\frac{1}{2}\int\limits_{B}
       \left[2V  e^{2u}+ \langle x,\nabla \log{V} \rangle Ve^{2u}\right] dx.
\end{eqnarray*}
However \eqref{eq:Diri} readily implies $\left(\,\int\limits_{\pa B}u_\nu \, d \s \right)^2=\rho^2$,
while we clearly have
$\int\limits_{B} V e^{2u}dx=\rho$. At this point an explicit calculation yields
$$
   \frac{1}{4\pi}\rho^2\leq -\frac{1}{2}\int\limits_{\pa B}(x,\nu)V e^{2u} d \s +\frac{1}{2}2 \rho+\frac{1}{2}2\al\rho<(1+\al)\rho,
   $$
and the conclusion follows. Observe that the sharpness of the strict inequality is due to the negative sign of the first term on the right in the first inequality which in fact vanishes along the well known radial and explicit solutions blowing up at the origin.
\end{pf}

\

\subsection{A non existence result on $\sd$ with two antipodal singularities}
We generalize an argument in \cite{tardcds} to obtain a non existence result for \eqref{eq:e-1}
in case of the sphere with two antipodal singularities.

\begin{pro}\label{p:nonexsph}
Let $(\Sigma, g) = (S^2, g_0)$, where $g_0$ is the standard round metric, let
$h\equiv 1$, $m=2$, let  $0<\al_1< \al_2<+\i$
be the weights of two antipodal singularities $\{p_1,p_2\}\subset \sd$ which we assume to coincide with the south and north pole respectively $p_1=S$, $p_2=N$.

Then a necessary condition for the solvability of \eqref{eq:e-1} is that
\begin{equation}\label{eq:nec}
  \hbox{ either } \quad 0<\rho<4\pi(1+\al_1), \qquad \qquad \hbox{ or } \quad
  4\pi(1+\al_2)<\rho<+\i.
\end{equation}
\end{pro}

\begin{pf}
We will work in isothermal coordinates induced by the stereographic projection $\Pi:\sd\mapsto\R^2$ satisfying
$\Pi(S)=0$. The local expression of the unique solution of \eqref{eq:Green} with $p=p_1=S$ takes the form
$$
G_S(\Pi^{-1}(z))=\frac{1}{4\pi}\log{\left(\frac{1+|z|^2}{2|z|^2}\right)}-\frac{1}{2\pi}\log{\left(\frac{e}{2}\right)}.
$$
In particular the local expression of the Laplace-Beltrami operator for the standard metric on
$\sd$ is
$$
\Delta_{g}=e^{-v_0}\Delta,
$$
where $\Delta$ is the standard Laplace operator in cartesian coordinates in $\R^2$ and $v_0$ satisfies
$$
v_0(z)=2\log{\left(\frac{2}{1+|z|^2}\right)},\qquad -\Delta v_0= 2 e^{v_0}\quad\mbox{in}\quad\R^2.
$$
Using these facts, and setting $\rho=2\pi\beta$, it is straightforward to check that $u$ solves \eqref{eq:e-1}
if and only if
\begin{eqnarray*}
  v(z) & =& 2 \, u(\Pi^{-1}(z)) + G_S(\Pi^{-1}(z)) + \frac{\beta-\al_2}{2}v_0(z)+
  2\al_1\log{\left(\frac{e}{2}\right)} \\ & + &
  (2+\al_2+\al_1-\beta)\log{2}+\log{(2\rho)} -\log{\int_{\sd} e^{2u} dV_{g_0}},
\end{eqnarray*}
satisfies
\begin{equation}\label{eq:vlocal}
\begin{cases}
 -\Delta v = K(z) e^{v} & \mbox{ in } \R^2; \\
  \int\limits_{\R^2}K(z)e^{v}dx = 4\pi\beta, &
\end{cases} \qquad \quad \hbox{ where } \quad
K(z)=\frac{|z|^{2\al_1}}{(1+|z|^2)^{2+\al_1+\al_2-\beta}}.
\end{equation}
Therefore we see that the results in \cite{clquant} can be applied to $v$ to conclude
$$
\int\limits_{\R^2} \langle z, \nabla K(z) \rangle e^{v} dx =4\pi\beta(\beta-2),
$$
so that, by using the integral constraint in \eqref{eq:vlocal}, an explicit evaluation shows
\begin{equation}\label{eq:Poho1}
2(2+\al_1+\al_2-\beta)\int\limits_{\R^2} \frac{|z|^2}{1+|z|^2}K(z) e^{v} dx
=4\pi\beta(2(1+\al_1)-\beta).
\end{equation}
Next, by writing $\frac{|z|^2}{1+|z|^2}=1-\frac{1}{1+|z|^2}$, and by using \eqref{eq:Poho1},
we obtain the independent constraint
\begin{equation}\label{eq:Poho2}
2(2+\al_1+\al_2-\beta)\int\limits_{\R^2} \frac{1}{1+|z|^2}K(z) e^{v}dx =4\pi\beta(2(1+\al_2)-\beta).
\end{equation}
By using \eqref{eq:Poho1} and \eqref{eq:Poho2} together and by discussing the cases
$2+\al_1+\al_2-\beta\lesseqqgtr 0$ it is readily seen that if $\al_1<\al_2$ then a necessary condition for the
solvability of \eqref{eq:vlocal} (and then of \eqref{eq:e-1}) is \eqref{eq:nec}.

Of course, by setting $\al_1=0$ we recover the non existence result obtained in \cite{tardcds} for the case where only
one singularity is considered.
\end{pf}

\begin{rem}
(a) Concerning the case $\al_1=0$ it has been already observed in \cite{tardcds} that in particular one
obtains in this way another proof of the non existence of conformal metrics with constant
Gaussian curvature on $\sd$ with one conical singularity, see \cite{Troy0} and the more recent
paper \cite{barjga}.
Indeed we obtain another proof of the non existence of conformal metrics with constant
Gaussian curvature on $\sd$ with two conical singularities of different orders $\al_1\neq\al_2$ which
corresponds to the case $2+\al_1+\al_2-\beta=0$. In fact in this situation
we see that \eqref{eq:Poho1} and \eqref{eq:Poho2} together imply $\al_1=\al_2$, and in this case solutions are
classified explicitly, see \cite{Troy0} and \cite{pratar}. The non existence results
for $2+\al_1+\al_2-\beta=0$ are associated with a well known problem, see \cite{Troy0}, corresponding to the best
pinching constants for these singular surfaces. The case with negative singularities has been
recently solved in \cite{barjga} while, at least to our knowledge, the positive case is  still open.

(b) We expect that existence should hold in some cases for which $\rho > 4 \pi \min_i \{1 + \al_i\}$.
For example we speculate that our method, with some extra work, could be adapted to the following situation: $m = 2$, $4 k \pi \leq \al_1, \al_2 < 4 (k+1) \pi$ for some $k \in \N$ and $\rho \in \left( 4 \pi \max \{1+ \al_1, 1+ \al_2\}, 4 (k+1) \pi \right)$.
\end{rem}


\begin{thebibliography}{99}


\bibitem{aubin} Aubin T., {\em Meilleures constantes dans le theoreme d'inclusion de Sobolev
et un theoreme de Fredholm non lineaire pour la transformation
conforme de la courbure scalaire}, J. Funct. Anal. 32 (1979),
148-174.




\bibitem{barjga} Bartolucci D., {\em On the best pinching constant of conformal metrics on $\sd$ with
one and two conical singularities}, Jour. Geom. Analysis (2012) to appear.

\bibitem{bclt} Bartolucci D.,  Chen C.C.,  Lin C.S.,  Tarantello G., {\em
Profile of blow-up solutions to mean field equations with singular data}, Comm. Part. Diff. Eq. 29 (2004), 1241-1265.

\bibitem{bardem} Bartolucci D., De Marchis F., {\em On the Ambjorn-Olesen electroweak condensates},
Jour. Math. Phys., to appear.

\bibitem{bdm} Bartolucci D., De Marchis F., Malchiodi A., {\em Supercritical
conformal metrics on surfaces with conical singularities},
Int. Math. Res. Not. (2011), Vol. 2011(24), 5625-5643.

\bibitem{barlin} Bartolucci D., Lin C.S., {\em Uniqueness
results for mean field equations with singular data},
Comm. Part. Diff. Eq.  34(7-9)  (2009), 676-702.

\bibitem{barlin-2} Bartolucci D., Lin C.S., {\em Sharp existence results for mean field equations with singular data},
Jour. Diff. Eq. 252(7) (2012), pp. 4115-4137.

\bibitem{barlintar} Bartolucci D., Lin C.S., Tarantello G., {\em Uniqueness and
symmetry results for solutions of a mean field equation on $\sd$ via a new bubbling
phenomenon}, Comm. Pure Appl. Math. 64(12) (2011),  1677-1730.

\bibitem{barmon-1} Bartolucci D., Montefusco E., {\em On the shape of blow-up solutions to a mean field
equation}, Nonlinearity 19 (2006), 611-631.

\bibitem{barmon} Bartolucci D., Montefusco E., {\em Blow up analysis,
existence and qualitative properties of solutions for the two
dimensional Emden-Fowler equation with singular potential},
M$^{2}$.A.S. 30(18) (2007), 2309--2327.


\bibitem{btcmp02} Bartolucci D., Tarantello G., {\em
Liouville type equations with singular data and their application
to periodic multivortices for the electroweak theory}, Comm. Math.
Phys. 229 (2002), 3-47.


\bibitem{bredon} Bredon G. E., Topology and Geometry. Grad. Texts in Math.
139, A. M. S., Providence, RI, 1997.

%

\bibitem{bls}  Brezis H., Li Y.Y., Shafrir I., {\em A $\sup+\inf$
inequality for some nonlinear elliptic equations involving
exponential nonlinearities}, J. Funct. Anal. 115 (1993), 344-358.

\bibitem{bm} Brezis H., Merle F., {\em Uniform estimates and blow-up
behavior for solutions of $-\Delta u =V(x) e\sp u$ in two
dimensions} Comm. Part. Diff. Eq. 16(8-9) (1991),
1223-1253.


\bibitem{clmp2} Caglioti E., Lions P.L., Marchioro C., Pulvirenti M.,
{\em A special class of stationary flows for two dimensional Euler
equations: a statistical mechanics description. II}, Comm. Math.
Phys.  174 (1995), 229-260.


\bibitem{carl} Carlotto A., {\em
On the solvability of singular Liouville equations on compact surfaces on arbitrary genus}, Trans. A.M.S. to appear.


\bibitem{cama} Carlotto A., Malchiodi A., {\em
Weighted barycentric sets and singular Liouville equations on compact surfaces},
J. Funct. Anal. 262(2) (2012), 409-450

\bibitem{cygc} Chang S.Y.A., Yang, P. C.,
{\em Conformal deformation of metrics on $S^2$}, J. Diff. Geom. 27 (1988), 259-296.

\bibitem{chakie} Chanillo S., Kiessling M., {\em Rotational symmetry of
solutions of some nonlinear problems in statistical mechanics and
in geometry}, Comm. Math. Phys. 160 (1994), 217-238.

\bibitem{chenwx}  Chen W.X., {\em A Trudinger inequality on surfaces with conical
singularities}, Proc. Amer. Math. Soc. 108 (1990), 821-832.

\bibitem{cl} Chen W., Li C., {\em Prescribing Gaussian curvatures on
surfaces with conical singularities}, J. Geom. Anal. 1(4) (1991),
359-372.

\bibitem{cl1} Chen C.C., Lin C.S., {\em Sharp estimates for solutions
of multi-bubbles in compact Riemann surfaces}, Comm. Pure Appl.
Math. 55(6) (2002), 728-771.

\bibitem{cl2} Chen C.C., Lin C.S., {\em Topological degree for a mean
field equation on Riemann surfaces}, Comm. Pure Appl. Math. 56(12)
(2003), 1667-1727.

\bibitem{cl4} Chen C.C., Lin C.S., Wang G., {\em Concentration phenomena of two-vortex solutions in a
Chern-Simons model}, Ann. Scuola Norm. Sup. Pisa Vol. III (2004), 367-397.

\bibitem{cl25} Chen C.C., Lin C.S., {\em Mean field equations of liouville type with
singular data: sharper estimates}, Discr. Cont. Dyn. Syt. 28(3) (2010), 1237-1272.

\bibitem{cl3} Chen C.C., Lin C.S., {\em A degree counting formulas for singular Liouville-type equation
and its application to multi vortices in electroweak theory}, in preparation.

\bibitem{clquant} Chen W., Li C., {\em Qualitative properties of solutions of
some nonlinear elliptic equations in $R^2$}, Duke Math. J. 71(2) (1993), 427-439.




\bibitem{cia} Cianchi A., {\em Moser-Trudinger trace inequalities}, Adv. Math. 217 
(2008), no. 5, 2005-2044. 


\bibitem{dpem} Del Pino M., Esposito P., Musso M., {\em Two-dimensional Euler
flows with concentrated vorticities}, Trans. Am. Math. Soc. 362(12) (2010), 6381–6395.


\bibitem{dem} De Marchis F., {\em Multiplicity result for a scalar field equation on
compact surfaces}, Comm. Part. Diff. Eq. 33(10-12) (2008), 2208-2224.


\bibitem{dem2} De Marchis F.,{\em Generic multiplicity for a scalar field equation on compact surfaces},
J. Funct. Anal. (259) (2010), 2165-2192.

\bibitem{djlw} Ding W., Jost J., Li J., Wang G., {\em Existence results
for mean field equations}, Ann. Inst. Henri Poincar\'e, Anal. Non
Lin\'eaire 16(5) (1999), 653-666.

\bibitem{djlw2} Ding W., Jost J., Li J., Wang G., {\em An analysis of the two-vortex case in the Chern-
Simons-Higgs model}, Calc. Var. P.D.E. 7 (1998), 87–97.

\bibitem{dj} Djadli Z., {\em Existence result for the mean field problem
on Riemann surfaces of all genuses}, Comm. Contemp. Math.  10(2) (2008), 205-220.

\bibitem{dm} Djadli Z., Malchiodi A., {\em Existence of conformal metrics with
constant $Q$-curvature}, Ann. of Math.  168(3)  (2008), 813-858.


\bibitem{det} Dolbeault J., Esteban M.J., Tarantello G., {\em The role of Onofri type inequalities in the symmetry properties
of extremals for Caffarelli-Kohn-Nirenberg inequalities, in two space dimensions}.
Ann. Sc. Norm. Super. Pisa Cl. Sci. Vol. VII  (2008),  313-341.




\bibitem{dunne} Dunne G., Self-dual Chern-Simons Theories, Lecture Notes in
Physics (1995).

\bibitem{esp} Esposito P., {\em Blowup solutions for a Liouville equation with singular
data}, SIAM J. Math. Anal. 36 (2005), no. 4, 1310-1345.

\bibitem{Fon} Fontana L., {\em Sharp borderline Sobolev inequalities on compact Riemannian
Manifolds}, Comm. Math. Elv.  68 (1993), 415-454.


\bibitem{fm} Fontana L., Morpurgo C., {\em Adams inequalities on measure spaces}, 
Adv. Math. 226 (2011), no. 6, 5066-5119. 




\bibitem{hkp} Hong J., Kim Y., Pac P.Y., {\em Multivortex Solutions of the
Abelian Chern-Simons Theory}, Phys. Rev. Lett. 64 (1990),
2230-2233.


\bibitem{jawe} Jackiw R., Weinberg E.J., {\em Selfdual Chern-Simons vortices},
Phys. Rev. Lett. 64 (1990), 2234-2237.


\bibitem{kk} Kallel S., Karoui R., {\em Symmetric joins and weighted barycenters},
Adv. Nonlinear Stud. 11(1) (2011), 117-143.


\bibitem{lai} Lai C.H. (ed.), Selected Papers on Gauge Theory of Weak and Electromagnetic
Interactions, World Scientific, Singapore, 1981.


\bibitem{li}  Li Y.Y., {\em Harnack type inequality: The method of moving planes},
Commun. Math. Phys. 200(2) (1999), 421-444.


\bibitem{lisha} Li Y.Y., Shafrir I., {\em Blow-up analysis for solutions of $- \D u =
V e^u$ in dimension two}, Indiana Univ. Math. J. 43(4) (1994),
1255-1270.

\bibitem{linweiye} Lin C.S., Wei J.C., Ye D., {\em Classification and nondegeneracy of SU(n + 1) Toda system with
singular sources}, Inventiones Math. (2012), to appear.




\bibitem{linwang} Lin C.S., Wang C.L., {\em Elliptic functions, Green functions
and the mean field equations on tori}, Ann of Math.
172(2) (2010), 911-954.


\bibitem{lucia} Lucia M., {\em A deformation lemma with an application to
a mean field equation}.  Topol. Methods Nonlinear Anal.  30(1)
(2007), 113-138.


\bibitem{luzh} Lucia M., Zhang L., {\em A priori estimates and uniqueness for
some mean field equations},  J. Diff. Eq.  217(1)  (2005), 154-178.


\bibitem{maldan} Malchiodi A., {\em Topological methods for an elliptic
equations with exponential nonlinearities}, Discrete Contin. Dyn.
Syst. 21(1) (2008), 277-294.


\bibitem{mald} Malchiodi A., {\em Morse theory and a scalar field equation
on compact surfaces}, Adv. Diff. Eq. 13(11-12)  (2008), 1109-1129.

\bibitem{maru} Malchiodi A., Ruiz D.,
{\em New improved Moser-Trudinger inequalities and singular Liouville equations on compact surfaces}, G.A.F.A., 21-5 (2011), 1196-1217.

\bibitem{maru2} Malchiodi A., Ruiz D., {\em A variational analysis of the toda system on compact surfaces},
Comm. Pure Appl. Math., to appear.



\bibitem{mc} McOwen R., {\em Conformal metrics in $\R^2$ with prescribed Gaussian curvature
and positive total curvature}, Indiana Univ. Math. J., 34 (1985), 97-104.


\bibitem{moser} Moser J., {\em A sharp form of an inequality by
N.Trudinger}, Indiana Univ. Math. J. 20 (1971), 1077-1091.


\bibitem{mo} Moser J., {\em On a nonlinear problem in differential geometry},
Dynamical Systems (M. Peixoto ed.), Academic Press, New York,
1973, 273-280.


\bibitem{ns} Nagasaki K., Suzuki T., {\em Asymptotic analysis for
two-dimensional elliptic eigenvalue problems with exponentially
dominated nonlinearities}.  Asymptotic Anal.  3(2)  (1990), 173-188.


%
%
%


\bibitem{pratar} Prajapat J., Tarantello G., {\em On a class of elliptic
problems in $R^{2}$: Symmetry and Uniqueness results},  Proc. Roy. Soc. Edinburgh 131A (2001), 967-985.

\bibitem{noltarcmp} Nolasco M., Tarantello G., {\em Vortex condensates for the SU(3) Chern-Simons theory},
Comm. Math. Phys., 213 (2000), 599-639.

\bibitem{osuz} Ohtsuka H., Suzuki T., {\em Blow-up analysis for SU(3) Toda system},
Jour. Diff. Eq. 232(2) (2007), 419-440.

\bibitem{sy2} Spruck J., Yang Y., {\em On Multivortices in the Electroweak Theory I:Existence of Periodic Solutions},
Comm. Math. Phys. 144 (1992), 1-16.


\bibitem{str} Struwe M., {\em The existence of surfaces of constant mean
curvature with free boundaries}, Acta Math. 160(1-2) (1988), 19-64.


\bibitem{st}  Struwe M., Tarantello G., {\em On multivortex solutions in
Chern-Simons gauge theory}, Boll. Unione Mat. Ital., Sez. B,
Artic. Ric. Mat. 8(1) (1998), 109-121.



\bibitem{tjfa05} Tarantello G., {\em A quantization property for blow up
solutions of singular Liouville-type equations}, J. Func. Anal. 219 (2005), 368-399.

\bibitem{tind05} Tarantello G., {\em  An Harnack inequality for
Liouville-type equations with singular sources},  Indiana Univ. Math J. 54(2) (2005), 599-615.

\bibitem{tar} Tarantello G., {\em Self-Dual Gauge Field Vortices:
An Analytical Approach}, PNLDE 72, Birkh\"auser Boston, Inc.,
Boston, MA, 2007.


\bibitem{tardcds} Tarantello G., {\em
Analytical, geometrical and topological aspects of a class of mean field equations on
surfaces}, Discrete Contin. Dyn.
Syst. 28(3) (2010), 931-973.

\bibitem{Troy0} Troyanov M., {\em Metrics of constant curvature on a sphere with two
conical singularities}, Proc. Third Int. Symp. on Diff. Geom. (Peniscola 1988),
Lect. Notes in Math. {\bf 1410} Springer-Verlag, 296--308.

\bibitem{tro} Troyanov M., {\em Prescribing Curvature on Compact Surfaces
with Conical Singularities}, Trans. Am. Math. Soc. 324(2) (1991), 793-821.

\bibitem{turya} Tur A., Yanovsky V., {\em Point vortices with a rational necklace: New exact stationary solutions
of the two-dimensional Euler equation}, Phys. Fluids 16(8) (2004), 2877-2885.


\bibitem{yang} Yang Y., Solitons in Field Theory and Nonlinear Analysis,
Springer Monographs in Mathematics, Springer, New York, 2001.

\bibitem{Zeidler} Zeidler E., {\em Nonlinear Functional Analysis and its Applications},
Springer-Verlag, New York, 1986.

\bibitem{luzh2} Zhang L., {\em Asymptotic behavior of blowup
solutions for elliptic equations with exponential nonlinearity and
singular data}, Comm. Contemp. Math. 11(3) (2009), 395-411.




\end{thebibliography}
\end{document}